\documentclass[preprint,12pt]{elsarticle}

\usepackage{amssymb}
\usepackage{amsmath}
\usepackage{graphicx}
\usepackage{subcaption}
\usepackage{float}
\usepackage{hyperref} 
\usepackage{color} 
\usepackage{scalerel,xparse}
\usepackage[normalem]{ulem}
\usepackage{amsthm}
\usepackage{ragged2e}

\begin{document}

\begin{frontmatter}



\title{A hat function based positive contact discontinuity capturing Boltzmann scheme}
\tnotetext[label1]{corresponding author}
\author[label2]{Balwinder Singh\corref{label1}}
\ead{balwinders@iisc.ac.in}
\author[label2]{S.V. Raghurama Rao}
\ead{raghu@iisc.ac.in}
\affiliation[label2]{organization={Dept. of Aerospace Engineering},
            addressline={Indian Institute of Science},
            city={Bengaluru},
            postcode={560012},
            country={India}}
            
\begin{abstract}
The suitable moments of the Maxwell velocity distribution function recover the conserved variable and flux vectors of the Euler equations of gas dynamics. However, the Maxwellian has infinite support in velocity space, which poses challenges for the numerical analysis of kinetic schemes based on it. Although Maxwellian is a natural choice, it is not the only distribution function that preserves these key moments. In his seminal paper, Perthame \cite{Perthame} introduced a compactly supported \textit{hat function} as an alternative. This paper presents the formulation and analysis of a contact discontinuity capturing Boltzmann scheme based on peculiar velocity, constructed using the hat function. Owing to the compact support of the hat function in velocity space, the numerical analysis of the proposed scheme is considerably simplified, allowing for a direct proof of its positivity. The scheme, extended to second-order accuracy without compromising the positivity-preserving property of its first-order counterpart, is tested on several benchmark inviscid flow problems in one- and two-dimensions.

\end{abstract}




\end{frontmatter}


\section{Introduction}
\label{sec1}

Extreme flow conditions, such as shock formation, are ubiquitous in high-speed flows. Thus, the prediction of high-speed flows requires robust, accurate, and efficient numerical methods. Boltzmann schemes represent an elegant alternative to traditional Riemann solvers. Such schemes are based on the fact that the Euler equations of gas dynamics are suitable moments of the Boltzmann equation of the kinetic theory of gases, whenever the velocity distribution function is a Maxwellian. Some of the prominent kinetic schemes for Euler equations are described in \cite{Chu}, \cite{Sanders_Prendergast}, \cite{Pullin}, \cite{Reitz}, \cite{SMD_NASA_TP_2613}, \cite{SMD_AIAA_1986}, \cite{SMD_NASA_TP_2583}, \cite{Kaniel}, \cite{Perthame}, \cite{Prendergast_Xu}, \cite{JCMandal_SMD}, \cite{SVRRao_SMD}.
\par
Of particular interest to us is the peculiar velocity based upwind (PVU) Boltzmann scheme \cite{SVRRao_SMD, SVRRao_SMD_91FM11, SVRRao_PVU_ICNMFD-94, SMD_Balak_SVRRao}. In the original work, separating fluid and peculiar velocities (and introducing upwinding based on each) is a desired feature leading to convection-pressure splitting at the macroscopic level. The original scheme is overly diffusive and cannot recognize any stationary discontinuities. Building on this framework, Singh and Raghurama Rao \cite{Singh_Rao_2026} recognize that an upwind discretization of the convective flux and a central discretization of the pressure flux can lead to exact resolution of an isolated stationary contact discontinuity.  However, this approach proved unstable in the presence of strong nonlinear discontinuities. To overcome this, the authors designed the peculiar velocity based upwind numerical flux such that the numerical diffusion corresponding to it switches off across stationary contact discontinuities and remains low but non-zero otherwise.

\par
The Maxwell velocity distribution function is infinitely supported in the velocity space, which causes difficulties in performing the numerical analyses on the Maxwellian-based kinetic schemes. Although the natural choice for the velocity distribution function is Maxwellian, it is not a unique distribution function whose moments lead to the conserved variable and flux vectors of the Euler equations of gas dynamics. Perthame \cite{Perthame} proposed a hat function with bounded support in the velocity space for the distribution function that preserves these key moments. This work presents the formulation and analysis of the modified PVU Boltzmann scheme based on the hat function, which notably captures stationary linear discontinuities. Owing to the compact support in velocity space of the hat function, the numerical analysis of the modified PVU Boltzmann scheme is considerably simplified, allowing for a direct proof of its positivity , while still producing numerical results equivalent to those obtained with the Maxwellian.
\par
This paper is organized as follows. Section \ref{sec 2} presents a brief overview of the kinetic model for polyatomic gases. In Section \ref{sec 3}, we study the moments of the Maxwellian distribution function and note that only four of the six moments are independent. We then derive the compactly supported hat function that satisfies these four moment constraints. We note that it is an alternative derivation to the one provided in \cite{Perthame}. In Section \ref{sec 4}, we present the numerical formulation of the one-dimensional peculiar velocity based upwind Boltzmann scheme using the hat function and establish its positivity-preserving property. In Section \ref{sec 5}, we present the numerical formulation and analysis of the contact discontinuity resolving PVU Boltzmann scheme using the hat function. In subsection \ref{subsec 6.1}, we study the moments of the two-dimensional Maxwellian distribution function and note that only eight of the twelve moments are independent. In the next subsection \ref{subsec 6.2}, we construct a two-dimensional hat function that satisfies the eight independent moments. Section \ref{sec 7} extends the novel kinetic scheme to two dimensions. Section \ref{sec 8} presents a way to extend the proposed kinetic scheme to second-order accuracy in space and time. Section \ref{sec 9} tests the scheme on the benchmark one- and two-dimensional inviscid gas dynamics test cases. Finally, the paper is concluded with remarks.


\section{From kinetic theory to macroscopic fluid dynamics}
\label{sec 2}

It is well known from the kinetic theory of gases that the continuum equations of fluid dynamics naturally emerge in the asymptotic limit of the Boltzmann equation; see, for instance, \cite{Saint-Raymond_2009}. In this section, we present a brief overview of the Boltzmann equation and its moments, focusing in particular on obtaining the compressible Euler equations. 
\par
The kinetic theory of gases represents a flow by the molecular velocity distribution function $ f(\boldsymbol{x}, \boldsymbol{v}, t, I) $ of particles at location $ \boldsymbol{x} $, time $ t $, with velocity $ \boldsymbol{v} $. An internal energy variable $ I $ corresponding to non-translational degrees of freedom is introduced to force the given value of $ \gamma $ (ratio of specific heats) \cite{Morse}. The Boltzmann equation governs the evolution of the velocity distribution function in an eight-dimensional phase space, which reads
\begin{equation} \label{3D_Boltzmann_Equation}
    \frac{\partial f}{\partial t} + \boldsymbol{v} \cdot \nabla f = Q(f, f),
\end{equation}
where $ Q(f, f) $ denotes the collision term. The velocity distribution function and collision term satisfy the following properties:
\begin{subequations} \label{equation 2}
    \begin{equation} \label{U_moments}
        \int_{\boldsymbol{v}} \, D\boldsymbol{v} \int_I \, dI \, \boldsymbol{\psi} f(\boldsymbol{x}, \boldsymbol{v}, t, I) = \boldsymbol{U}(\boldsymbol{x}, t) \quad \forall f,
    \end{equation}
    \vspace{-0.35 cm}
    \begin{equation} \label{moment_of_collision_term}
        \int_{\boldsymbol{v}} \, D\boldsymbol{v} \int_I \, dI \, \boldsymbol{\psi} \, Q(f,f) = 0 \quad \forall f,
    \end{equation}
    \begin{equation} \label{moment_of_log_times_collision_term}
        \int_{\boldsymbol{v}} \, D\boldsymbol{v} \int_I \, dI \, \log{f} \, Q(f,f) \leq 0 \quad \forall f,
    \end{equation}
\end{subequations}
where $ \boldsymbol{U} $ denotes the vector of conserved variables and $ \boldsymbol{\psi} $ denotes the moment vector defined as
\begin{equation} \label{definition_of_moment_vector}
    \boldsymbol{\psi} (\boldsymbol{v}, I) \equiv \begin{bmatrix}
        1 & \boldsymbol{v} & I + \dfrac{1}{2} |\boldsymbol{v}|^2
    \end{bmatrix}^T.
\end{equation}
We note that a detailed structure of the collision term is not required in view of (\ref{moment_of_collision_term}). Integrating the Boltzmann equation (\ref{3D_Boltzmann_Equation}) with respect to $ \boldsymbol{v}$ and $I$ leads to
\begin{equation} \label{moment_of_the_Boltzmann_equation}
    \int_{\boldsymbol{v}} \, D\boldsymbol{v} \int_I \, dI \, \boldsymbol{\psi} (\boldsymbol{v}, I) \biggl\{ {\frac{\partial f}{\partial t}} + \boldsymbol{v} \cdot \nabla f \biggl\} = 0.
\end{equation}
At equilibrium, $ Q(f, f) = 0 $, and the unique solution of this equation is the Maxwell distribution function, which reads
\begin{equation} \label{3D_Maxwell_distribution_function}
    f^{eq} = \dfrac{\rho}{I_0} \biggl( \dfrac{\beta}{\pi} \biggl)^{(D_f/2)} \exp{\biggl( -\beta \bigl( \boldsymbol{v} - \boldsymbol{u} \bigl)^2 - \dfrac{I}{I_0} \biggl)},
\end{equation}
where $ \rho $ denotes density, $ \boldsymbol{u} $ denotes bulk fluid velocity, $ D_f $ denotes translational degrees of freedom, $ \beta = 1/(2 R T) $, $ T $ denotes temperature, and
\begin{equation} \label{definition_of_I_0}
    I_0 = \frac{(2+D_f)-\gamma D_f}{4(\gamma-1)\beta}.
\end{equation}
The moments of $ \boldsymbol{v} f^{eq} $ recover the inviscid flux tensor $ \boldsymbol{\mathcal{G}} $, and thus, the moments of the Boltzmann equation in the limit $ f \rightarrow f^{eq} $ recover the Euler equations of gas dynamics, \textit{i.e.},
\begin{equation} \label{Boltzmann_to_Euler}
    \int_{\boldsymbol{v}} \, D\boldsymbol{v} \int_{I} \, dI \, \boldsymbol{\psi} \, \biggl( {\frac{\partial f}{\partial t}} + \boldsymbol{v}\cdot\nabla f = 0, f = f^{eq} \biggl) \Leftrightarrow \frac{\partial \boldsymbol{U}}{\partial t} + \nabla \cdot \boldsymbol{\mathcal{G}} = 0,
\end{equation}
where $ \boldsymbol{U} $ and $ \boldsymbol{\mathcal{G}} $ are defined as
\begin{subequations}
    \begin{align}
        \boldsymbol{U} &= \begin{bmatrix}
            \rho & \rho \boldsymbol{u} & \rho E
        \end{bmatrix}, \\
        \boldsymbol{\mathcal{G}} &= \begin{bmatrix}
            \rho \boldsymbol{u} & \rho \boldsymbol{u} \otimes \boldsymbol{u} + p & \rho \boldsymbol{u} H
        \end{bmatrix}.
    \end{align}
\end{subequations}
Here
\begin{equation} \label{definition_of_E_e}
    E = \dfrac{p}{\rho (\gamma - 1)} + \dfrac{| \boldsymbol{u} |^2}{2}, \; H = E + \frac{p}{\rho},
\end{equation} 
and $ p $ denotes the thermodynamic pressure. We refer to \cite{SMD_NASA_TP_2613}, \cite{SMD_NASA_TP_2583}, \cite{Vincenti_Kruger} for more details on the connection between the macroscopic equations and the Boltzmann equation of the kinetic theory of gases.
\subsection{Peculiar velocity based Boltzmann equation}
The peculiar velocity (also known as {\em thermal velocity} or {\em relative velocity}) $ \textbf{c} $ is defined to be the relative velocity of a particle with respect to the bulk fluid velocity $ \textbf{u} $, \textit{i.e.},
\begin{equation} \label{peculiar_velocity_definition}
    \boldsymbol{c} \equiv \boldsymbol{v} - \boldsymbol{u}.
\end{equation}
A peculiar velocity based Boltzmann equation is obtained by replacing the molecular velocity with a sum of peculiar and bulk fluid velocities in (\ref{3D_Boltzmann_Equation}), which in $ k $ spatial dimensions reads
\begin{equation} \label{peculiar_velocity_based_Boltzmann_equation}
    \frac{\partial f}{\partial t} + \nabla \cdot \left( \boldsymbol{u} f \right) + \nabla \cdot \left( \boldsymbol{c} f \right) = Q(f, f).
\end{equation}
We take note of some interesting properties of the peculiar velocity based Boltzmann equation. We redefine the moment vector in terms of the peculiar and bulk fluid velocity as
\begin{equation} \label{moment_vector_in_terms_of_peculiar_velocity}
    \boldsymbol{\psi} (\boldsymbol{u} (\boldsymbol{x}, t), \boldsymbol{c} (\boldsymbol{x}, t), I) = \begin{bmatrix}
        1 & \left( \boldsymbol{u} + \boldsymbol{c} \right) & I + \dfrac{1}{2} | \boldsymbol{u} + \boldsymbol{c} |^2
    \end{bmatrix}^T,
\end{equation}
and take moments of (\ref{peculiar_velocity_based_Boltzmann_equation}) in the limit $ f \to f^{eq} $ to obtain a system of flux-split Euler equations, which reads
\begin{equation} \label{Zha-Bilgen_type_split_Euler_equations}
    \frac{\partial \textbf{U}}{\partial t} + \nabla \cdot \boldsymbol{\mathcal{G}}^{(c)} + \nabla \cdot \boldsymbol{\mathcal{G}}^{(p)} = 0,
\end{equation}
where $ \boldsymbol{\mathcal{G}}^{(c)} $ and $ \boldsymbol{\mathcal{G}}^{(p)} $ represent the $ k $ dimensional tensor of advective and pressure fluxes respectively. The $ i^{th} $ components of these tensors are defined by the following moment relations:
\begin{subequations} \label{definitions_of_convective_pressure_fluxes}
    \begin{align}
        \textbf{G}_i^{(c)} &= \int_{I} \, dI \int_{\textbf{c}} \, d\textbf{c} \, \boldsymbol{\psi} \, u_i \, f^{eq} = \begin{bmatrix}
        \rho u_i & \rho u_i u_j & \rho u_i E 
        \end{bmatrix}^T \quad \forall i, \\
        \textbf{G}_i^{(p)} &= \int_{I} \, dI \int_{\textbf{c}} \, d\textbf{c} \, \boldsymbol{\psi} \, c_i \, f^{eq} = \begin{bmatrix}
        0 & p \delta_{i j} & p u_i
        \end{bmatrix}^T \quad \forall i.
    \end{align}
\end{subequations}
The peculiar velocity based Boltzmann equation separates the fluid velocity based uni-directional (along the tangent to a streamline) information propagation from the random velocity based multi-directional information propagation (corresponding to the jiggling of the molecules in random directions). We refer to \cite{SVRRao_SMD, SMD_Balak_SVRRao, Singh_Rao_2026} for the description of all the features of the peculiar velocity based Boltzmann equation, corresponding upwind method and the physics of the convection-pressure splitting.


\section{A compactly supported distribution function in one dimension}
\label{sec 3}

\subsection{Moments of the one-dimensional Maxwell distribution function}
\label{subsec3.1}

We first consider the moment relations corresponding to the conserved variable vector. The zeroth moment of the Maxwell distribution function recovers the first component of the conserved variable vector, \textit{i.e.},
\begin{equation*}
    U_1 = \int_{0}^{\infty} \, dI \int_{-\infty}^{\infty} \, dc \, {f^{eq}} .
\end{equation*} 
Integrating the above equation with respect to $ I $ leads to
\begin{equation} \label{Zeroth_F_integral}
    \int_{-\infty}^{\infty} \, dc \, \bar{f}^{eq} = \rho, 
\end{equation}
where $ \bar{f}^{eq} $ denotes the integrated Maxwell distribution function with respect to $ I $, which is written as
\begin{equation} \label{Definition_of_Maxwellian_integrated_wrt_I}
    \bar{f}^{eq} = {\rho} \sqrt{\dfrac{\beta}{\pi}} \exp{\left( -\beta \, c^2 \right)},
\end{equation}
and $ c $ denotes the one-dimensional peculiar velocity defined as
\begin{equation}
    c = v - u.
\end{equation}
Note that $ dv = dc $. The first moment of the Maxwell distribution function recovers the second component of the conserved variable vector, \textit{i.e.},
\begin{equation*}
    U_2 = \int_0^{\infty} \, dI \int_{-\infty}^{\infty} \, dc \, \left( u + c \right) {f^{eq}} = \rho u.
\end{equation*}
Again integrating the above equation w.r.t $ I $ leads to
\begin{equation} \label{First_F_integral}
    \int_{-\infty}^{\infty} \, dc \, c \, \bar{f}^{eq} = 0.
\end{equation}
The second moment of the Maxwell distribution function reads
\begin{equation*}
    U_3= \int_0^{\infty} \, dI \int_{-\infty}^{\infty} \, dc \, \Bigl(I + \dfrac{1}{2} u^2 + \dfrac{1}{2} c^2 + u c \Bigl){f^{eq}} = \rho E.
\end{equation*}
Integrating with respect to $ I $ and rearranging the terms give
\begin{equation} \label{Second_F_integral}
    \int_{-\infty}^{\infty} \, dc \, c^2 \, \bar{f}^{eq} = p.
\end{equation}
We now consider the moments of the Maxwell distribution function corresponding to the flux vector. We note that the moments recovering the first and second components of the flux vector yield the same relations as those in equations (\ref{First_F_integral}) and (\ref{Second_F_integral}).
Finally, we consider the moment relation for the energy flux $ \boldsymbol{G}_1^{(3)} $, \textit{i.e.},
\begin{equation*}
    \boldsymbol{G}_1^{(3)}(\boldsymbol{U}) = \int_0^{\infty} \, dI \int_{-\infty}^{\infty} \, dc \, \Bigl(I + {\dfrac{1}{2}} u^2 + \dfrac{1}{2} c^2 + u c \Bigl) (u + c) {f^{eq}} = pu + \rho u E.
\end{equation*}
Integrating the above equation with respect to $ I $ and rearranging the terms gives
\begin{equation} \label{Third_F_integral}
    \int_{-\infty}^{\infty} \, dc \, c^3 \, \bar{f}^{eq} = 0.
\end{equation}
The study of all six moments corresponding to the conserved variable and flux vectors of the Euler equations of gas dynamics reveals that there are only four independent moment relations in one dimension. Any distribution function $ \bar{f}^{eq} $ that satisfies these four integrals will, in turn, recover Euler equations for $ \gamma = 5/3 $. We list the four independent moment relations below for future reference:
\begin{subequations} \label{Reduced_F_moments}
    \begin{align}
        \int_{-\infty}^{\infty} \, dc \, \bar{f}^{eq} &= \rho, \label{Reduced_F_moment_1} \\
        \int_{-\infty}^{\infty} \, dc \, c \bar{f}^{eq} &= 0, \label{Reduced_F_moment_2} \\
        \int_{-\infty}^{\infty} \, dc \, c^2 \bar{f}^{eq} &= p, \label{Reduced_F_moment_3} \\
        \int_{-\infty}^{\infty} \, dc \, c^3 \bar{f}^{eq} &= 0. \label{Reduced_F_moment_4}
    \end{align}
\end{subequations}

\subsection{A compactly supported hat function in one dimension}
\label{subsec3.2}
\begin{figure}[h!]
   \centering
   \includegraphics[scale=0.6, trim={0 0 0 0}, clip]{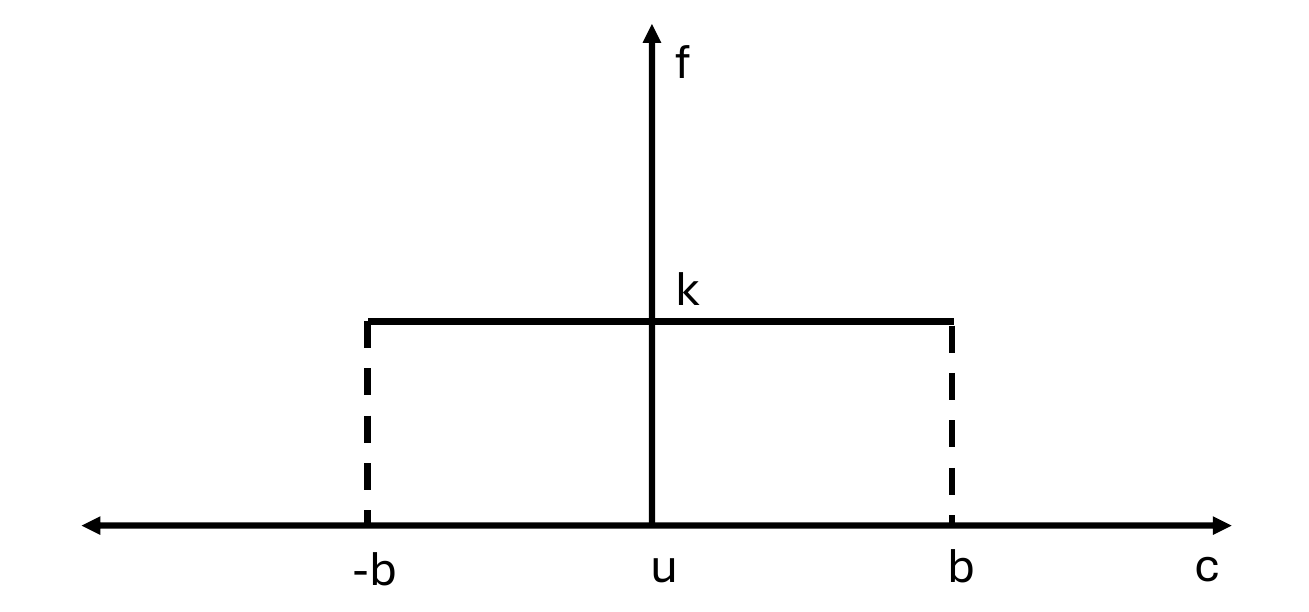}
   \caption{A one-dimensional hat function}
   \label{fig: one-dimensional hat function}
\end{figure}

In this subsection, we construct a compactly supported distribution function, with its parameters uniquely determined by enforcing the four independent moment relations as listed in equation (\ref{Reduced_F_moments}).
We replace the Maxwell distribution function with a compactly supported hat function as shown in Figure \ref{fig: one-dimensional hat function}.
\noindent
The compactly supported distribution function is assumed to be of the form
\begin{equation} \label{Proposed_1D_distribution_function}
    \bar{f}^{eq} = 
    \begin{cases}
        0 &\quad \text{if} \quad c < -b \\[10pt]
        \rho \, k &\quad \text{if} \quad -b \leq c \leq b \\[10pt]
        0 &\quad \text{if} \quad c > b \\
    \end{cases}
\end{equation}
The maximum value $ k $ and the half-width $ b $ of the hat function are obtained by enforcing it to satisfy the four independent moment relations. We first substitute equation (\ref{Proposed_1D_distribution_function}) into equation (\ref{Reduced_F_moment_1}) to obtain
\begin{flalign} \label{First_relation}
    2 \, k \, b = 1,
\end{flalign}
and then in equation (\ref{Reduced_F_moment_3}) to obtain
\begin{flalign} \label{Definition_of_half-width}
    b = \sqrt{3 R T}.
\end{flalign}
Substituting the value of half-width from equation (\ref{Definition_of_half-width}) into equation (\ref{First_relation}) leads to
\begin{flalign}
    k = \dfrac{1}{2 \, \sqrt{3 R T}}.
\end{flalign}
Remark: The compactly supported hat function defined in (\ref{Proposed_1D_distribution_function}) satisfies the moment relations (\ref{Reduced_F_moment_2}) and (\ref{Reduced_F_moment_4}) identically.
\par
\noindent
Thus, the compactly supported hat function is defined as
\begin{equation} \label{Definition_of_new_1D_linear_distribution_function}
    \bar{f}^{eq} = 
    \begin{cases}
        \quad 0 &\quad \text{if} \quad c < -\sqrt{3 R T} \\[10pt]
        \dfrac{\rho}{2 \sqrt{3 R T}} &\quad \text{if} \quad -\sqrt{3 R T} \leq c \leq \sqrt{3 R T} \\[10pt]
        \quad 0 &\quad \text{if} \quad c > \sqrt{3 R T} \\
    \end{cases}
\end{equation}

\noindent
\textbf{Note:} This derivation is an alternative to the one originally proposed by Perthame in his seminal paper \cite{Perthame}.

\subsection{Equivalent form of the Euler equations of gas dynamics in 1D}
\label{subsec3.3}

The moments of the compactly supported equilibrium distribution function defined in (\ref{Definition_of_new_1D_linear_distribution_function}) recover the conserved variable and the flux vectors of the Euler equations for monatomic gases, \textit{i.e.},
\begin{subequations}
    \begin{align}
        \int_c \, dc \, \Tilde{\boldsymbol{\psi}} \, \bar{f}^{eq} &= \begin{bmatrix}
            \rho & \rho u & \dfrac{1}{2} \rho u^2 + \dfrac{1}{2} p
        \end{bmatrix}^T, \\
        \int_c \, dc \, \Tilde{\boldsymbol{\psi}} \, \left(u + c \right) \, \bar{f}^{eq} &= \begin{bmatrix}
            \rho u & p + \rho u^2 & \dfrac{1}{2} \rho u^3 + \dfrac{1}{2} p u
        \end{bmatrix}^T,
    \end{align}
\end{subequations}
where $ \Tilde{\boldsymbol{\psi}} $ is defined as
\begin{equation}
    \Tilde{\boldsymbol{\psi}} = \begin{bmatrix}
        1 & \left( u + c \right) & \dfrac{1}{2} |u + c|^2
    \end{bmatrix}^T.
\end{equation}
One way to take into account the internal degrees of freedom of molecules was introduced in Section (\ref{sec 2}). However, the resulting distribution function is unbounded in the internal energy space $ I $. To address this, Perthame later proposed an alternative approach in \cite{Perthame1999}, introducing two compactly supported distribution functions, $ \bar{f}^{eq} $ and $ \bar{g}^{eq} $, such that
\begin{equation}
    \int_c \, dc \, \dfrac{|u + c|^2}{2} \bar{f}^{eq} + \int_c \, dc \, \bar{g}^{eq} = \dfrac{1}{2} \rho u^2 + \dfrac{p}{\left( \gamma - 1 \right)}.
\end{equation}
A simple choice for $ \bar{g}^{eq} $ is
\begin{equation} \label{second_distribution_function}
    \bar{g}^{eq} = I_0 \bar{f}^{eq},
\end{equation}
where $ I_0 $ is defined in (\ref{definition_of_I_0}) with $ D_f = 1 $. Adopting this selection allows the conserved quantities and the associated flux vectors of the Euler equations governing gas dynamics to be written in the following form:
\begin{subequations}
    \begin{align}
        \boldsymbol{U} &= \int_c \, dc \, \Tilde{\boldsymbol{\psi}}_1 \, \bar{f}^{eq} + \int_c \, dc \, \Tilde{\boldsymbol{\psi}}_2 \, \bar{g}^{eq}, \label{1D_definition_of_conserved_variable} \\
        \boldsymbol{G}_1 &= \int_c \, dc \, \Tilde{\boldsymbol{\psi}}_1 \, \left( u + c \right) \, \bar{f}^{eq} + \int_c \, dc \, \Tilde{\boldsymbol{\psi}}_2 \, \left( u + c \right) \, \bar{g}^{eq},
    \end{align}
\end{subequations}
where $ \Tilde{\boldsymbol{\psi}}_1 $ and $ \Tilde{\boldsymbol{\psi}}_2 $ are defined as
\begin{equation}
    \Tilde{\boldsymbol{\psi}}_1 = \begin{bmatrix}
        1 & \left( u + c \right) & \dfrac{1}{2} |u + c|^2
    \end{bmatrix}^T, \; \Tilde{\boldsymbol{\psi}}_2 = \begin{bmatrix}
        0 & 0 & 1
    \end{bmatrix}^T.
\end{equation}
The flow is now described by two velocity distribution functions, $ f(x, u, c, t) $ and $ g(x, u, c, t) $, at spatial location $ x $ and time $ t $, with macroscopic velocity $ u $ and peculiar velocity $ c $. These functions evolve according to the following Boltzmann-like equations:
\begin{subequations} \label{Boltzmann-like_equations}
    \begin{align}
        \dfrac{\partial f}{\partial t} + \dfrac{\partial}{\partial x} \left( u f \right) + \dfrac{\partial}{\partial x} \left( c f \right) = Q_1, \\
        \dfrac{\partial g}{\partial t} + \dfrac{\partial}{\partial x} \left( u g \right) + \dfrac{\partial}{\partial x} \left( c g \right) = Q_2,
    \end{align}
\end{subequations}
where $ Q_1 $ and $ Q_2 $ are unknown functions satisfying for all $ t $, $ x $
\begin{subequations}
    \begin{align}
        \int_c \, dc \, Q_1 = 0, \\
        \int_c \, dc \, \left( u + c \right) \, Q_1 = 0, \\
        \int_c \, dc \, \left( \dfrac{|u + c|^2}{2} \, Q_1 + Q_2 \right) = 0.
    \end{align}
\end{subequations}
Thus, in the limits $ f \rightarrow \bar{f}^{eq} $ and $ g \rightarrow \bar{g}^{eq} $, the Euler equations of gas dynamics can be recovered as a linear combination of the two Boltzmann-like equations (\ref{Boltzmann-like_equations}), \textit{i.e.},
\begin{multline}
    \int_c \, dc \, \Tilde{\boldsymbol{\psi}}_1 \left( \dfrac{\partial f}{\partial t} + \dfrac{\partial}{\partial x} \left( u f \right) + \dfrac{\partial}{\partial x} \left( c f \right), \; f = \bar{f}^{eq} \right) \\
    + \int_c \, dc \, \Tilde{\boldsymbol{\psi}}_2 \left( \dfrac{\partial g}{\partial t} + \dfrac{\partial}{\partial x} \left( u g \right) + \dfrac{\partial}{\partial x} \left( c g \right), \; g = \bar{g}^{eq} \right) = 0.
\end{multline}


\section{The one-dimensional peculiar velocity based upwind Boltzmann scheme with a hat function}
\label{sec 4}

In this section, we present a one-dimensional, first-order, hat function-based peculiar velocity based upwind (PVU) Boltzmann scheme for the numerical solutions of the Euler equations of gas dynamics. A one-dimensional domain is discretized into a finite number of control volume cells, indexed as $ j $. Cell boundaries are labeled as $ x_{j \pm \frac{1}{2}} $. Initially, the bulk fluid velocity and the thermodynamic state variables are assumed to be piecewise constant in each finite volume cell, \textit{i.e.},
\begin{subequations} \label{constant_state_distribution}
    \begin{align}
        u(x,0) &= u_j \quad \forall x \in  (x_{j-\frac{1}{2}}, x_{j+\frac{1}{2}}), \\
        \rho(x,0) &= \rho_j \quad \forall x \in  (x_{j-\frac{1}{2}}, x_{j+\frac{1}{2}}), \\
        T(x,0) &= T_j \quad \forall x \in (x_{j-\frac{1}{2}}, x_{j+\frac{1}{2}}).
    \end{align}
\end{subequations}
Consequently, the equilibrium distribution functions defined in (\ref{Definition_of_new_1D_linear_distribution_function}) and (\ref{second_distribution_function}) are constructed inside each cell. The distribution functions in each finite volume cell are evolved by numerically solving the peculiar velocity based Boltzmann equations (\ref{Boltzmann-like_equations}) in two stages:
\par
\textbf{Convection step:} In this stage, the particles are allowed to stream freely from one cell to their neighbor cells. We start with the one-dimensional collision-less peculiar velocity based Boltzmann equations, which reads
\begin{subequations} \label{1D_collisionless_peculiar_velocity_based_Boltzmann_equation}
    \begin{align}
        \frac{\partial f}{\partial t} + \frac{\partial}{\partial x} (u f) + \frac{\partial}{\partial x} (c f) = 0, \\
        \frac{\partial g}{\partial t} + \frac{\partial}{\partial x} (u g) + \frac{\partial}{\partial x} (c g) = 0.
    \end{align}
\end{subequations}
We use a classical upwind scheme for the numerical discretization of (\ref{1D_collisionless_peculiar_velocity_based_Boltzmann_equation}), which reads
\begin{subequations} \label{discretized_PV_Boltzmann_equation}
    \begin{align}
        f_j^{(n + 1)-} - f_j^{(n)} + \sigma \left[ \hat{f}_{j + \frac{1}{2}}^{(n)} - \hat{f}_{j - \frac{1}{2}}^{(n)} \right] = 0, \label{discretized_PV_Boltzmann_equation_f} \\
        g_j^{(n + 1)-} - g_j^{(n)} + \sigma \left[ \hat{g}_{j + \frac{1}{2}}^{(n)} - \hat{g}_{j - \frac{1}{2}}^{(n)} \right] = 0, \label{discretized_PV_Boltzmann_equation_g}
    \end{align}
\end{subequations}
where $ \sigma = \Delta t/\Delta x $, $ \hat{f}_{j \pm \frac{1}{2}}^{(n)} $ and $ \hat{g}^{(n)}_{j \pm \frac{1}{2}} $ are the numerical fluxes (defined at the kinetic scale) at the cell interfaces, which reads
\begin{subequations}
    \begin{align}
        \hat{h}_{j + \frac{1}{2}} &= (u_j^+ + c_j^+) \Bar{h}^{eq}_j + (u_{j + 1}^- + c_{j + 1}^-) \Bar{h}^{eq}_{j + 1}, \\
        \hat{h}_{j - \frac{1}{2}} &= (u_{j - 1}^+ + c_{j - 1}^+) \Bar{h}^{eq}_{j - 1} + (u_{j}^- + c_{j}^-) \Bar{h}^{eq}_{j},
    \end{align}
\end{subequations}
for $ h = $ $ f $ or $ g $.
\par
\textbf{Collision step:} The streaming step drives the distribution functions $ f $ and $ g $ away from the local thermodynamic equilibrium.
In the collision stage, the effect of particle collisions is taken into account by instantaneously relaxing the post-streaming distribution functions towards a new equilibrium state. This relaxation is performed in such a way that the conserved quantities-mass, momentum, and total energy within each cell-remain unchanged. It is equivalent to solving the following ordinary differential equations in the Euler limit:
\begin{equation} \label{equation_equivalent_collision_step}
    \frac{d f}{d t} = Q_1, \; \frac{d g}{d t} = Q_2,
\end{equation}
\textit{i.e.},
\begin{subequations} \label{essence_of_collision_step}
    \begin{align}
        \int_{c} \, dc \, \Tilde{\boldsymbol{\psi}}_1 \, \Bar{f}_j^{eq, (n + 1)} = \int_{c} \, dc \, \Tilde{\boldsymbol{\psi}}_1 \, f_j^{(n + 1)-}, \\
        \int_{c} \, dc \, \Tilde{\boldsymbol{\psi}}_2 \, \Bar{g}_j^{eq, (n + 1)} = \int_{c} \, dc \, \Tilde{\boldsymbol{\psi}}_2 \, g_j^{(n + 1)-}.
    \end{align}
\end{subequations}
Using Equations (\ref{discretized_PV_Boltzmann_equation}) and (\ref{essence_of_collision_step}) in (\ref{1D_definition_of_conserved_variable}) leads to a first-order finite volume scheme for the Euler equations of gas dynamics, which reads
\begin{equation} \label{1D_finite_volume_update_formula}
    \bar{\boldsymbol{U}}_j^{n + 1} = \bar{\boldsymbol{U}}_j^{n} - \sigma \left( \hat{\textbf{f}}_{j + \frac{1}{2}} - \hat{\textbf{f}}_{j - \frac{1}{2}} \right),
\end{equation}
where $ \bar{\boldsymbol{U}} $ denotes the cell-averaged value of the conserved variable vector in each cell, $ \hat{\textbf{f}}_{j \pm \frac{1}{2}} $ denotes the interface numerical flux, which is defined as the sum of advective and pressure numerical fluxes. The numerical advective and pressure fluxes are defined as follows:
\begin{multline} \label{definition_advective_numerical_flux}
    \hat{\textbf{f}}_{j + \frac{1}{2}}^{(a)} = \int_{-\infty}^{\infty} \, dc \, \Tilde{\boldsymbol{\psi}}_1 \, \left( u_j^+ \, \Bar{f}_j^{eq} + u_{j + 1}^- \, \Bar{f}_{j + 1}^{eq} \right) \\
    + \int_{-\infty}^{\infty} \, dc \, \Tilde{\boldsymbol{\psi}}_2 \, \left( u_j^+ \, \Bar{g}_j^{eq} + u_{j + 1}^- \, \Bar{g}_{j + 1}^{eq} \right),
\end{multline}
\begin{multline} \label{definition_pressure_numerical_flux}
    \hat{\textbf{f}}_{j + \frac{1}{2}}^{(p)} = \int_{-\infty}^{\infty} \, dc \, \Tilde{\boldsymbol{\psi}}_1 \, \left( c_j^+ \, \Bar{f}_j^{eq} + c_{j + 1}^- \, \Bar{f}_{j + 1}^{eq} \right) \\
    + \int_{-\infty}^{\infty} \, dc \, \Tilde{\boldsymbol{\psi}}_2 \, \left( c_j^+ \, \Bar{g}_j^{eq} + c_{j + 1}^- \, \Bar{g}_{j + 1}^{eq} \right),
\end{multline}
where 
\begin{equation} \label{CIR_splitting_of_peculiar_fluid_velocities}
    u^\pm = \frac{u \pm |u|}{2}, \, c^\pm = \frac{c \pm |c|}{2}.
\end{equation}
Remark: Notice that the collision step introduces the temporal discontinuity in the velocity distribution function only and not the conserved variable vector, \textit{i.e.}, 
\begin{equation}
    \Bar{f}^{(n+1)-} \neq \Bar{f}^{eq, (n+1)}, \, \Bar{\textbf{U}}^{(n+1)-} = \Bar{\textbf{U}}^{(n+1)}.
\end{equation}

\subsection{Positivity analysis in 1D}
\label{subsec4.1}
\noindent
\textbf{Theorem 1.} Assume $ \Bar{f}^{eq}_j = 0 $ for $ |c_j| \geq |c_j|_{max} $, then under the CFL-like condition
\begin{equation} \label{CFL-like_condition}
    \max_{j} \{ |u_j| + |c_j|_{max} \} \Delta t \leq \Delta x,
\end{equation}
the peculiar velocity based upwind Boltzmann scheme (\ref{1D_finite_volume_update_formula})-(\ref{CIR_splitting_of_peculiar_fluid_velocities}) is positive, \textit{i.e.}, $ \rho \geq 0, p \geq 0 $.

\vspace{1\baselineskip}
\noindent
Proof: We rewrite (\ref{discretized_PV_Boltzmann_equation_f}) using (\ref{definition_advective_numerical_flux}) and (\ref{definition_pressure_numerical_flux}) in finite difference form as follows:
\begin{equation}
    f_j^{(n+1)-} = \Bar{f}_j^{eq, (n)} (1 - \sigma (|u_j| + |c_j|)) + \sigma (h_{j-1}^+ + h_{j+1}^-),
\end{equation}
with
\begin{subequations} \label{Defintions_of_h+_and_h-}
    \begin{align}
        h_{j-1}^+ &= \left( u^+ \, \Bar{f}^{eq, (n)} \right)_{j - 1} + \left( c^+ \, \Bar{f}^{eq, (n)} \right)_{j - 1}, \\
        h_{j+1}^- &= -\left( u^- \, \Bar{f}^{eq, (n)} \right)_{j + 1} - \left( c^- \, \Bar{f}^{eq, (n)} \right)_{j + 1}.
    \end{align}
\end{subequations}
We note that
\begin{subequations}
    \begin{align}
        h_{j-1}^+ &\geq 0, \\
        h_{j+1}^- &\geq 0.
    \end{align}
\end{subequations}
We only have to consider the peculiar velocities for which $ \Bar{f}^{eq}_j > 0 $; these are given by
\begin{equation} \label{peculiar_velocities_for_which_f_is_positive}
    |c_j| \leq |c_j|_{max}.
\end{equation}
We add $ |u_j| $ to both sides of (\ref{peculiar_velocities_for_which_f_is_positive}) and note that
\begin{equation}
    \sigma (|u_j| + |c_j|) \leq \sigma (\max_{j} \{ |u_j| + |c_j|_{max} \}) \leq 1.
\end{equation}
Hence, we have $ f_j^{(n + 1)-} \geq 0 $ and thus
\begin{subequations}
    \begin{align}
        \rho_j^{(n+1)} &= \int_{c} \, dc \, f_j^{(n+1)-} \geq 0, \\
        p_j^{(n+1)} &= \int_{c} \, dc \, c_j^2 \, f_j^{(n+1)-} \geq 0,
    \end{align}
\end{subequations}
and Theorem 1. is proved.

\subsection{Inability of the PVU Scheme to capture stationary contact discontinuities}
\label{subsec 4.2}

Consider a stationary contact discontinuity. The left and right states satisfy
\begin{equation} \label{contact discontinuity}
    u_L = u_R = 0, \; \text{and} \;\; p_L = p_R = p_0, \; \text{with} \;\; \rho_L \neq \rho_R.
\end{equation}
The exact solution corresponding to this initial condition is
\begin{equation}
    \boldsymbol{U} (x, t) = \boldsymbol{U} (x, 0) \quad \forall t.
\end{equation}
Therefore, a stationary contact discontinuity is an exact equilibrium of the Euler equations. A numerical method that claims to capture stationary contact discontinuities exactly should admit the same equilibrium at the discrete level. A sufficient condition for discrete equilibrium is that every numerical interface flux assumes the same value, \textit{i.e.,}
\begin{equation}
    \hat{\textbf{f}}_{j + \frac{1}{2}} = \begin{bmatrix}
        0 & p_0 & 0
    \end{bmatrix}^T \quad \forall j.
\end{equation}

In the PVU scheme, the numerical flux can be written as
\begin{equation}
    \hat{\textbf{f}}_{j + \frac{1}{2}} = \dfrac{1}{2} \left( \textbf{G}_j + \textbf{G}_{j + 1} \right) - \textbf{d}_{j + \frac{1}{2}}^{(a)} - \textbf{d}_{j + \frac{1}{2}}^{(p)}.
\end{equation}
For the stationary contact discontinuity, the average flux reduces to
\begin{equation} \label{average flux at stationary contact}
    \dfrac{1}{2} \left( \textbf{G}_j + \textbf{G}_{j + 1} \right) = \begin{bmatrix}
        0 & p_0 & 0
    \end{bmatrix}^T
\end{equation}
and the advective diffusive flux $ \boldsymbol{d}_{j + \frac{1}{2}}^{(a)} $ vanishes. Therefore, the stationary contact discontinuity is smeared only because
\begin{equation}
    \textbf{d}_{j + \frac{1}{2}}^{(p)} \neq 0.
\end{equation}


\section{A contact discontinuity capturing PVU Boltzmann scheme}
\label{sec 5}
The peculiar velocity based upwind Boltzmann scheme based on the hat function is positive and robust, but unable to recognize stationary linear discontinuities. In addition to the diffusive prediction of a stationary contact discontinuity, the PVU Boltzmann scheme also suffers from spurious pressure and velocity oscillations near isolated contact discontinuities. However, the Boltzmann equation based on peculiar velocity (\ref{peculiar_velocity_based_Boltzmann_equation}) provides an elegant framework for treating peculiar velocity-based flux independently of the bulk fluid velocity-based flux. In this section, we attempt to design a peculiar velocity-based numerical diffusion that vanishes as the Mach number goes to zero and remains nonzero otherwise, thus improving the scheme's performance across the linear discontinuities without compromising the robustness of the original scheme.

\subsection{Numerical formulation}
\label{subsec 5.1}

We choose a general form of the numerical flux for the advective and pressure parts which reads
\begin{multline} \label{definition_advective_numerical_flux_PVU_plus}
    \hat{\textbf{f}}_{j + \frac{1}{2}}^{(a)} = \int_c \, dc \, \Tilde{\boldsymbol{\psi}}_1 \left( \frac{1}{2} (u \bar{f}^{eq})_j + \frac{1}{2} (u \bar{f}^{eq})_{j+1} - \frac{1}{2} \Delta_{j + \frac{1}{2}} \left( \phi^{(u)}_{j + \frac{1}{2}} |u| \bar{f}^{eq} \right) \right) \\
    + \int_c \, dc \, \Tilde{\boldsymbol{\psi}}_2 \left( \frac{1}{2} (u \bar{g}^{eq})_j + \frac{1}{2} (u \bar{g}^{eq})_{j+1} - \frac{1}{2} \Delta_{j + \frac{1}{2}} \left( \phi^{(u)}_{j + \frac{1}{2}} |u| \bar{g}^{eq} \right) \right),
\end{multline}
\begin{multline} \label{definition_pressure_numerical_flux_PVU_plus}
    \hat{\textbf{f}}_{j + \frac{1}{2}}^{(p)} = \int_c \, dc \, \Tilde{\boldsymbol{\psi}}_1 \left( \frac{1}{2} (c \bar{f}^{eq})_j + \frac{1}{2} (c \bar{f}^{eq})_{j+1} - \frac{1}{2} \Delta_{j + \frac{1}{2}} \left( \phi^{(c)}_{j + \frac{1}{2}} |c| \bar{f}^{eq} \right) \right) \\
    + \int_c \, dc \, \Tilde{\boldsymbol{\psi}}_2 \left( \frac{1}{2} (c \bar{g}^{eq})_j + \frac{1}{2} (c \bar{g}^{eq})_{j+1} - \frac{1}{2} \Delta_{j + \frac{1}{2}} \left( \phi^{(c)}_{j + \frac{1}{2}} |c| \bar{g}^{eq} \right) \right),
\end{multline}
with
\begin{equation} \label{definition_of_delta}
    \Delta_{j + \frac{1}{2}} (\cdot) = (\cdot)_{j + 1} - (\cdot)_{j}.
\end{equation}
In addition to the above, $ \phi's $ must also obey the dynamical consistency condition, \textit{i.e.},
\begin{equation} \label{dynamical_consistency_of_phi's}
    \textbf{U}_j = \textbf{U}_{j \pm 1} \implies \left( \phi^{(u), (c)}_{j \pm \frac{1}{2}} \right)_j = \left( \phi^{(u), (c)}_{j \pm \frac{1}{2}} \right)_{j \pm 1}.
\end{equation}

\subsection{Positivity analysis}
\label{subsec 5.2}

In this subsection, we study the positivity-preserving property of the modified PVU Boltzmann scheme while keeping $ \phi $'s general.

\vspace{1\baselineskip}
\textbf{Lemma 2.}
In addition to the dynamical consistency condition (\ref{dynamical_consistency_of_phi's}), suppose that $ \phi^{(c)} $ satisfies
\begin{subequations} \label{first_set_of_constraints_on_phi_2}
\begin{align}
\text{(a)} \quad & \phi^{(c)} \geq 0, \\
\text{(b)} \quad & \phi^{(c)}(-c) = \phi^{(c)}(c)
\quad \text{or} \quad
\phi^{(c)} \neq \phi^{(c)}(c).
\end{align}
\end{subequations}
Then the following inequalities hold:
\begin{subequations} \label{Definitions_of_J1_and_J2}
    \begin{align}
        J_1 &= \int_{-\infty}^{\infty} \, dc \, \frac{(c + |c| \phi^{(c)})}{2} \bar{f}^{eq} \geq 0, \label{Definition_of_J1} \\
        J_2 &= -\int_{-\infty}^{\infty} \, dc \, \frac{(c - |c| \phi^{(c)})}{2} \bar{f}^{eq} \geq 0. \label{Definition_of_J2}
    \end{align}
\end{subequations}

\vspace{1\baselineskip}
Proof: We rewrite (\ref{Definition_of_J1}) as
\begin{equation} \label{splitting_J1_into_J1+_J1-}
    J_1 = \int_{-\infty}^{0} \, dc \frac{(c + |c| \phi^{(c)})}{2} \Bar{f}^{eq} + \int_{0}^{\infty} \, dc \frac{(c + |c| \phi^{(c)})}{2} \Bar{f}^{eq}
\end{equation}
We substitute $ c $ with $ -c $ in the first term of (\ref{splitting_J1_into_J1+_J1-}) to obtain
\begin{equation} \label{substituting_-c_into_J1}
     \Bar{J}_1 = \int_{0}^{\infty} \, dc \frac{(-c + |c| \phi^{(c)})}{2} \Bar{f}^{eq} + \int_{0}^{\infty} \, dc \frac{(c + |c| \phi^{(c)})}{2} \Bar{f}^{eq}.
\end{equation}
We then combine the two integrals to obtain
\begin{equation}
    J_1 = \int_{0}^{\infty} \, dc \, |c| \, \phi^{(c)} \, \Bar{f}^{eq} \geq 0.
\end{equation}
We now consider the positivity of (\ref{Definition_of_J2}) and rewrite it as
\begin{equation} \label{splitting_J2_into_J2+_J2-}
    J_2 = -\int_{-\infty}^{0} dc \frac{(c - |c| \phi^{(c)})}{2} \Bar{f}^{eq} - \int_{0}^{\infty} dc \frac{(c - |c| \phi^{(c)})}{2} \Bar{f}^{eq}.
\end{equation}
We substitute $ c $ with $ -c $ in the second term of (\ref{splitting_J2_into_J2+_J2-}) to obtain
\begin{equation} \label{substituting_-c_into_J2}
    J_2 = -\int_{-\infty}^{0} dc \frac{(c - |c| \phi^{(c)})}{2} \Bar{f}^{eq} - \int_{-\infty}^{0} dc \frac{(-c - |c| \phi^{(c)})}{2} \Bar{f}^{eq}.
\end{equation}
Finally, the simplification of (\ref{substituting_-c_into_J2}) leads to
\begin{equation}
    J_2 = \int_{-\infty}^{0} \, dc \, |c| \, \phi^{(c)} \, \Bar{f}^{eq} \geq 0.
\end{equation}
Hence, Lemma 2. is proved.

\vspace{1\baselineskip}
\textbf{Lemma 3.}
With $ \phi^{(c)} $ satisfying the constraints in (\ref{first_set_of_constraints_on_phi_2}), the following inequalities can be established:
\begin{subequations} \label{Definitions_of_J3_J4}
    \begin{align}
        J_3 &= \int_{-\infty}^{\infty} \, dc \, c^2 \frac{(c + |c| \phi^{(c)})}{2} \bar{f}^{eq} \geq 0, \label{Definition_of_J3} \\
        J_4 &= -\int_{-\infty}^{\infty} \, dc \, c^2 \frac{(c - |c| \phi^{(c)})}{2} \bar{f}^{eq} \geq 0. \label{Definition_of_J4}
    \end{align}
\end{subequations}

\vspace{1\baselineskip}
Proof: We rewrite (\ref{Definition_of_J3}) as
\begin{equation} \label{splitting_J3_into_J3+_J3-}
    J_3 = \int_{-\infty}^{0} \, dc \, c^2 \frac{(c + |c| \phi^{(c)})}{2} \Bar{f}^{eq} + \int_{0}^{\infty} \, dc \, c^2 \frac{(c + |c| \phi^{(c)})}{2} \Bar{f}^{eq}.
\end{equation}
We replace $ c $ with $ -c $ in the first term of (\ref{splitting_J3_into_J3+_J3-}) to get
\begin{equation} \label{substituting_-c_into_J3}
    J_3 = \int_{0}^{\infty} \, dc \, c^2 \frac{(-c + |c| \phi^{(c)})}{2} \Bar{f}^{eq} + \int_{0}^{\infty} \, dc \, c^2 \frac{(c + |c| \phi^{(c)})}{2} \Bar{f}^{eq}.
\end{equation}
We simplify (\ref{substituting_-c_into_J3}) and note that
\begin{equation}
    J_3 = \int_{0}^{\infty} \, dc \, c^2 \, |c| \, \phi^{(c)} \, \Bar{f}^{eq} \geq 0.
\end{equation}
We now rewrite (\ref{Definition_of_J4}) as
\begin{equation} \label{splitting_J4_into_J4+_J4-}
    J_4 = -\int_{-\infty}^{0} \, dc \, c^2 \frac{(c - |c| \phi^{(c)})}{2} \Bar{f}^{eq} - \int_{0}^{\infty} \, dc \, c^2 \frac{(c - |c| \phi^{(c)})}{2} \Bar{f}^{eq}.
\end{equation}
We substitute $ c $ with $ -c $ in the second term of (\ref{splitting_J4_into_J4+_J4-}) to obtain
\begin{equation} \label{substituting_-c_into_J4}
    J_4 = -\int_{-\infty}^{0} \, dc \, c^2 \frac{(c - |c| \phi^{(c)})}{2} \Bar{f}^{eq} - \int_{-\infty}^{0} \, dc \, c^2 \frac{(-c - |c| \phi^{(c)})}{2} \Bar{f}^{eq}.
\end{equation}
We finally collect both integrals together to get
\begin{equation}
    J_4 = \int_{-\infty}^{0} \, dc \, c^2 \, |c| \, \phi^{(c)} \, \Bar{f}^{eq} \geq 0.
\end{equation}
Hence, Lemma 3. is proved.

\vspace{1\baselineskip}
\textbf{Lemma 4.} 
In addition to (\ref{dynamical_consistency_of_phi's}), assume that $ \phi^{(u)} $ satisfies
\begin{subequations}
\begin{align}
\text{(a)} \quad & \phi^{(u)} > 0, \\
\text{(b)} \quad & \phi^{(u)} \neq \phi^{(u)}(c).
\end{align}
\end{subequations}
Then the following integrals satisfy
\begin{subequations} \label{Definitions_of_J5_and_J6}
    \begin{align}
        J_5 &= \int_{-\infty}^{\infty} \, dc \, \frac{(u + |u| \phi^{(u)})}{2} \bar{f}^{eq} \geq 0, \label{Definition_of_J5} \\
        J_6 &= -\int_{-\infty}^{\infty} \, dc \, \frac{(u - |u| \phi^{(u)})}{2} \bar{f}^{eq} \geq 0\label{Definition_of_J6}
    \end{align}
\end{subequations}
\textit{if and only if} $ \phi^{(u)} \geq 1 $.

\vspace{1\baselineskip}
Proof: Since $ \phi^{(u)} $ is not a function of the integration variable $ c $, therefore, the term $ u + | u | \phi^{(u)} $ can be factored out of the integral. Thus, Equation (\ref{Definition_of_J5}) becomes
\begin{equation} \label{re-writing J5}
    J_5 = \frac{(u + |u| \phi^{(u)})}{2} \int_{-\infty}^{\infty} \, dc \, \bar{f}^{eq}.
\end{equation}
Note that the integral
\begin{equation*}
    \int_{-\infty}^{\infty} \, dc \, \bar{f}^{eq} = \rho,
\end{equation*}
which is strictly positive. So, the integral (\ref{re-writing J5}) now reads
\begin{equation}
    J_5 = \rho \, \frac{(u + |u| \phi^{(u)})}{2}.
\end{equation}
Observe that for any real $ u $,
\begin{equation}
    u + | u | \phi^{(u)} = \begin{cases}
        u \, \left( 1 + \phi^{(u)} \right), \quad u \geq 0, \\
        u \, \left( 1 - \phi^{(u)} \right), \quad u < 0.
    \end{cases}
\end{equation}
\begin{itemize}
    \item If $ \phi^{(u)} \geq 1 $, then $ 1 - \phi^{(u)} \leq 0 $. For $ u < 0 $, we have $ u \left( 1 - \phi^{(u)} \right) \geq 0 $. For $ u \geq 0 $, $ u \left( 1 + \phi^{(u)} \right) \geq 0 $ trivially. Hence, $ u + | u | \phi^{(u)} \geq 0 $ for all $ u $, implying $ J_5 \geq 0 $.
    \item Conversely, if $ J_5 \geq 0 $ for all admissible distributions, then necessarily $ \phi^{(u)} \geq 1 $; otherwise, choosing $ u < 0 $ would make $ u \left( 1 - \phi^{(u)} \right) $ negative, yielding a negative integrand and therefore a negative $ J_5 $.
\end{itemize}
Thus, $ J_5 \geq 0 $ \textit{if and only if} $ \phi^{(u)} \geq 1 $. An analogous argument holds for $ J_6 $.

\vspace{1\baselineskip}
\textbf{Lemma 5.} 
Under the same assumptions on $ \phi^{(u)} $ as in Lemma 4, the following inequalities hold
\begin{subequations} \label{Definitions_of_J7_and_J8}
    \begin{align}
        J_7 &= \int_{-\infty}^{\infty} \, dc \, \frac{(u + |u| \phi^{(u)})}{2} \, c^2 \, \bar{f}^{eq} \geq 0, \label{Definition_of_J7} \\
        J_8 &= -\int_{-\infty}^{\infty} \, dc \, \frac{(u - |u| \phi^{(u)})}{2} \, c^2 \, \bar{f}^{eq} \geq 0\label{Definition_of_J8}
    \end{align}
\end{subequations}
\textit{if and only if} $ \phi^{(u)} \geq 1 $.

\vspace{1\baselineskip}
Proof: Since $ \phi^{(u)} $ is independent of the integration variable $ c $, $ u + | u | \phi^{(u)} $ can be factored out of the integral. Thus, Equation (\ref{Definition_of_J7}) becomes
\begin{equation} \label{re-writing J7}
    J_7 = \frac{(u + |u| \phi^{(u)})}{2} \int_{-\infty}^{\infty} \, dc \, c^2 \, \bar{f}^{eq}.
\end{equation}
Note that the integral
\begin{equation*}
    \int_{-\infty}^{\infty} \, dc \, c^2 \, \bar{f}^{eq} = p,
\end{equation*}
which is strictly positive. So, the integral (\ref{re-writing J7}) becomes
\begin{equation}
    J_7 = p \, \frac{(u + |u| \phi^{(u)})}{2}.
\end{equation}
As shown in Lemma 4,
\begin{equation}
    \dfrac{u + | u | \phi^{(u)}}{2} \geq 0 \quad \forall u,
\end{equation}
\textit{if and only if} $ \phi^{(u)} \geq 1 $. An identical argument applied to $ J_8 $ yields the same necessary and sufficient condition on $ \phi^{(u)} $.

\vspace{1\baselineskip}
\textbf{Note:} We further restrict $ \phi^{(c)} $ to the interval $ [0, 1] $. The condition $ \phi^{(c)} \leq 1 $ ensures reduced numerical diffusion associated with the peculiar velocity-based flux. Meanwhile, Lemmas 4 and 5 impose $ \phi^{(u)} \geq 1 $. Collectively, the two parameters must satisfy
\begin{subequations} \label{constraints_on_phi's}
    \begin{align}
        \phi^{(u)} &\geq 1, \label{constraint_on_phi1} \\ 
        0 \leq \phi^{(c)} &\leq 1. \label{constraint_on_phi2}
    \end{align}
\end{subequations}

\vspace{1\baselineskip}
\textbf{Theorem 6.} Assume that $ \bar{f}^{eq}_j = 0 $ for $ |c_j| > |c_j|_{max} $. Furthermore, assume that the diffusion control parameter $ \phi^{(u)} $, in addition to (\ref{constraint_on_phi1}) and (\ref{dynamical_consistency_of_phi's}), obeys
\begin{equation} \label{additional_constraint_on_phi1}
    \frac{1}{2} \frac{|u_j|}{\displaystyle \max_{\forall j} \{ |u_j| \}} \left( \phi^{(u)}_{j + \frac{1}{2}} + \phi^{(u)}_{j - \frac{1}{2}} \right)_j \leq m, \, m \in \mathbb{N}.
\end{equation}
Then, under the CFL-like condition,
\begin{equation} \label{CFL-like_condition_for_PVU_plus}
    \max_{j} \left\{ m \, |u_j| + |c_j|_{max} \right\} \Delta t \leq \Delta x,
\end{equation}
the modified PVU Boltzmann scheme defined by the fluxes (\ref{definition_advective_numerical_flux_PVU_plus}) and (\ref{definition_pressure_numerical_flux_PVU_plus}) is positive, \textit{i.e.,} $ \rho \geq 0, p \geq 0 $.

\vspace{1\baselineskip}
Proof: We have, from (\ref{1D_finite_volume_update_formula}), (\ref{definition_advective_numerical_flux_PVU_plus}), (\ref{definition_pressure_numerical_flux_PVU_plus}), (\ref{definition_of_delta}), and (\ref{constraints_on_phi's}), an update formula for the density, which reads
\begin{equation} \label{update_formula_for_density}
     \rho_j^{(n+1)} = \int_{c} \, dc \left( \bar{f}_j^{eq, (n)} (1 - \sigma (\zeta^{(1)}_j |u_j| + \zeta^{(2)}_j |c_j|)) + \sigma (h_{j-1}^+ + h_{j+1}^-) \right),
\end{equation}
with
\begin{subequations} \label{definitions_of_zeta's}
    \begin{align}
        \zeta^{(1)}_j &= \frac{1}{2} \left( \phi^{(u)}_{j + \frac{1}{2}} + \phi^{(u)}_{j - \frac{1}{2}} \right)_j, \\
        \zeta^{(2)}_j &= \frac{1}{2} \left( \phi^{(c)}_{j + \frac{1}{2}} + \phi^{(c)}_{j - \frac{1}{2}} \right)_j,
    \end{align}
\end{subequations}
and
\begin{subequations} \label{definition_of_h}
    \begin{align}
        h_{j-1}^+ &= \left( \Tilde{u}^+ \Bar{f}^{eq, (n)} \right)_{j - 1} + \left( \Tilde{c}^+ \Bar{f}^{eq, (n)} \right)_{j - 1}, \\
        h_{j+1}^- &= -\left( \Tilde{u}^- \Bar{f}^{eq, (n)} \right)_{j + 1} - \left( \Tilde{c}^- \Bar{f}^{eq, (n)} \right)_{j + 1}.
    \end{align}
\end{subequations}
We note that
\begin{subequations} \label{definitions_of_u_tilda_and_c_tilda}
    \begin{align}
        \Tilde{u}^+_{j - 1} &= \frac{u_{j - 1} + |u_{j - 1}| \left( \phi^{(u)}_{j - \frac{1}{2}} \right)_{j-1}}{2}, \label{positivity_of_u_tilda_plus} \\
        \Tilde{u}^-_{j + 1} &= \frac{u_{j + 1} - |u_{j + 1}| \left( \phi^{(u)}_{j + \frac{1}{2}} \right)_{j + 1}}{2}, \label{positivity_of_u_tilda_minus} \\
        \Tilde{c}_{j-1}^+ &= \frac{c_{j-1} + |c_{j-1}| \left( \phi^{(c)}_{j - \frac{1}{2}} \right)_{j - 1}}{2}, \\
        \Tilde{c}_{j+1}^- &= \frac{c_{j+1} - |c_{j+1}| \left( \phi^{(c)}_{j + \frac{1}{2}} \right)_{j + 1}}{2}.
    \end{align}
\end{subequations}
We then have, from Lemmas 2 and 4, that
\begin{subequations}
    \begin{align}
        \int_{c} \, dc \, h_{j-1}^+ &\geq 0, \\
        \int_{c} \, dc \, h_{j+1}^- &\geq 0.
    \end{align}
\end{subequations}
We also notice that
\begin{equation}
    \zeta^{(2)}_j \leq 1.
\end{equation}
We only have to consider the peculiar velocities for which $ \Bar{f}^{eq}_j > 0 $; these are given by
\begin{equation}
    |c_j| \leq |c_j|_{max}.
\end{equation}
We multiply $ |c_j| $ by $ \zeta^{(2)}_j $ and note that
\begin{equation} \label{equation_86}
    |c_j| \, \zeta^{(2)}_j \leq |c_j|_{max} \leq \max_{\forall j} |c_j|_{max}.
\end{equation}
We multiply $ |u_j| $ by $ \zeta^{(1)}_j $ and condition (\ref{additional_constraint_on_phi1}) dictates that
\begin{equation} \label{equation_87}
    |u_j| \, \zeta^{(1)}_j \leq m \max_{\forall j} |u_j|.
\end{equation}
Finally, we add equations (\ref{equation_86}) and (\ref{equation_87}) and by the CFL-like condition (\ref{CFL-like_condition_for_PVU_plus}), note that
\begin{equation} \label{equation_88}
    \sigma \left( |u_j| \zeta^{(1)}_j + |c_j| \zeta^{(2)}_j \right) \leq \sigma \left( m \max_{\forall j} |u_j| + \max_{\forall j} |c_j|_{max} \right) \leq 1.
\end{equation}
Hence, we have
\begin{equation}
    \int_{c} \, dc \, f_j^{(n+1)-} = \rho_j^{(n+1)} \geq 0.
\end{equation}
We now show that the modified PVU Boltzmann scheme based on peculiar velocity, defined by the fluxes (\ref{definition_advective_numerical_flux_PVU_plus}) and (\ref{definition_pressure_numerical_flux_PVU_plus}), ensures the positivity of pressure. We start with an update formula for the pressure, which reads
\begin{equation} \label{update_formula_for_pressure}
    p_j^{(n+1)} = \int_{c} \, dc \, c_j^2 \left( \bar{f}_j^{eq, (n)} (1 - \sigma (\zeta^{(1)}_j |u_j| + \zeta^{(2)}_j |c_j|)) + \sigma (h_{j-1}^+ + h_{j+1}^-) \right).
\end{equation}
We note, from Lemmas 3 and 5, that
\begin{subequations}
    \begin{align}
        \int_{c} \, dc \, c_j^2 \, h_{j-1}^+ &\geq 0, \\
        \int_{c} \, dc \, c_j^2 \, h_{j+1}^- &\geq 0.
    \end{align}
\end{subequations}
From (\ref{equation_88}), we note that
\begin{equation}
    \int_{c} \, dc \, c_j^2 \, \bar{f}_j^{eq, (n)} \left( 1 - \sigma (\zeta^{(1)}_j |u_j| + \zeta^{(2)}_j |c_j|) \right) \geq 0.
\end{equation}
Hence, $ p_j^{(n+1)} \geq 0 $ and Theorem 6. is proved.

\subsection{Contact-discontinuity constraints on the diffusion control parameter}
\label{subsubsec 5.3}

The numerical diffusion must satisfy
\begin{equation}
    \boldsymbol{d}_{j + \frac{1}{2}}^{(a)} = \boldsymbol{d}_{j + \frac{1}{2}}^{(p)} = 0   
\end{equation}
whenever
\begin{equation}
    u_L = u_R = 0
\end{equation}
and
\begin{equation}
    p_L = p_R.
\end{equation}

This condition becomes the additional defining principle of the diffusion design. Equivalently, the diffusion control parameter $ \phi^{(c)} $ must vanish on the stationary contact discontinuity manifold
\begin{equation}
    \mathcal{C} = \left\{ \left( \boldsymbol{U}_L, \boldsymbol{U}_R \right): u_L = u_R = 0, p_L = p_R \right\}
\end{equation}
Thus,
\begin{equation}
    \phi^{(c)} = 0 \quad \text{on} \quad \mathcal{C}.
\end{equation}

The diffusion control parameter $ \phi^{(u)} $ must, in addition, preserve the vanishing of the advective diffusive flux at stationary contact discontinuities.

\vspace{-1\baselineskip}
\subsection{Selection of diffusion control parameters}
\label{subsubsec 5.4}

The diffusion control parameter $ \phi^{(c)} $ should depend only on quantities that vanish for a stationary contact discontinuity. Natural candidates are
\begin{equation}
    \Delta p = p_R - p_L, \quad \Delta u = u_R - u_L, \quad \bar{u} = \dfrac{u_R + u_L}{2}.
\end{equation}
Thus, a general form for the diffusion coefficient is
\begin{equation}
    \phi^{(c)} = \Phi \left( |\Delta p|, |\Delta u|, |\bar{u}| \right),
\end{equation}
where
\begin{equation}
 \Phi(0, 0, 0) = 0.
\end{equation}

To investigate different functional forms of the contact-discontinuity diffusion control, three variants are considered. In all variants, the velocity-diffusion coefficient is defined as
\begin{equation} \label{phi_u}
    \left( \phi^{(u)}_{j + \frac{1}{2}} \right)_j = \dfrac{1}{|u_j|} \left( \max\left\{ |u_j|, |u_{j + 1}| \right\} + \left| \dfrac{\sqrt{\rho_j} u_j + \sqrt{\rho_{j + 1}} u_{j+1}}{\sqrt{\rho_j} + \sqrt{\rho_{j+1}}} \right| \right).
\end{equation}
The three variants differ only in the choice of $ \phi^{(c)} $, as described below.

\subsubsection{Variant 1}
\label{subsubsec:variant1}

For the first variant, the diffusion coefficient $ \phi^{(c)} $ is defined by
\begin{equation} \label{variant 1}
    \phi^{(c)}_{j + \frac{1}{2}} = \dfrac{S_c}{1 + S_c},
\end{equation}
where $ S_c $ is a dimensionless contact-discontinuity sensor given by
\begin{equation}
    S_c = \dfrac{|\Delta p|}{\bar{p}} + \dfrac{|\Delta u|}{\bar{a}},
\end{equation}
with
\begin{equation}
    \bar{p} = \dfrac{p_L + p_R}{2}, \quad \bar{a} = \dfrac{a_L + a_R}{2}.
\end{equation}

\subsubsection{Variant 2}
\label{subsubsec:variant2}

The second variant employs a hyperbolic-tangent form for $ \phi^{(c)} $:
\begin{equation} \label{variant 2}
    \phi^{(c)}_{j + \frac{1}{2}} = \tanh(S_c),
\end{equation}
where $ S_c $ is defined as in Variant 1.

\subsubsection{Variant 3}
\label{subsubsec:variant3}

As an alternative to the sensor-based forms above, the third variant employs an exponential form:
\begin{equation} \label{variant 3}
    \left(\phi^{(c)}_{j + \frac{1}{2}} \right)_j = \exp\left( -\dfrac{c_j^2}{u_j^2} \right).
\end{equation}

The one-dimensional numerical fluxes corresponding to $ \phi $'s in \eqref{phi_u} and \eqref{variant 3} are given in Appendix A.

\vspace{1\baselineskip}
\textbf{Theorem 7.} Let the left and right states satisfy the conditions of a stationary contact discontinuity given by \eqref{contact discontinuity}. Then the modified PVU Boltzmann scheme equipped with the diffusion control parameters defined in \eqref{phi_u}, \eqref{variant 1}, \eqref{variant 2}, and \eqref{variant 3} admits this discontinuity as a discrete equilibrium. In particular,
\begin{equation} \label{discrete equilibrium}
    \bar{\boldsymbol{U}}_j^{(n + 1)} = \bar{\boldsymbol{U}}_j^{(n)}
\end{equation}
holds for all $ j $ and all $ n \geq 0 $.

\vspace{1\baselineskip}
Proof: We have, from \eqref{1D_finite_volume_update_formula},
\begin{equation}
    \bar{\boldsymbol{U}}_j^{n + 1} = \bar{\boldsymbol{U}}_j^{n} - \sigma \left( \hat{\textbf{f}}_{j + \frac{1}{2}} - \hat{\textbf{f}}_{j - \frac{1}{2}} \right),
\end{equation}
where the interface numerical fluxes can be expressed as
\begin{subequations}
    \begin{align}
        \hat{\textbf{f}}_{j + \frac{1}{2}} &= \dfrac{1}{2} \left( \boldsymbol{G}_j + \boldsymbol{G}_{j + 1} \right) - \boldsymbol{d}_{j + \frac{1}{2}}^{(a)} - \boldsymbol{d}_{j + \frac{1}{2}}^{(p)}, \label{j + 1/2 flux} \\
        \hat{\textbf{f}}_{j - \frac{1}{2}} &= \dfrac{1}{2} \left( \boldsymbol{G}_{j - 1} + \boldsymbol{G}_j \right) - \boldsymbol{d}_{j - \frac{1}{2}}^{(a)} - \boldsymbol{d}_{j - \frac{1}{2}}^{(p)} \label{j - 1/2 flux}.
    \end{align}
\end{subequations}

From \eqref{definition_advective_numerical_flux_PVU_plus}, we note that
\begin{equation}
    d_{j + \frac{1}{2}}^{(a)} = \left( \phi_{j + \frac{1}{2}}^{(u)} \, | u | \, \boldsymbol{U} \right)_{j + 1} - \left( \phi_{j + \frac{1}{2}}^{(u)} \, | u | \, \boldsymbol{U} \right)_j
\end{equation}

From \eqref{phi_u}, we note that
\begin{equation}
    \lim_{u \rightarrow 0} \phi^{(u)} \, | u | = 0.
\end{equation}
Thus,
\begin{equation}
    \boldsymbol{d}_{j + \frac{1}{2}}^{(a)} = \boldsymbol{d}_{j - \frac{1}{2}}^{(a)} = 0.
\end{equation}
From \eqref{variant 1}, \eqref{variant 2}, and \eqref{variant 3}, we note that
\begin{equation}
    \lim_{u \rightarrow 0} \phi^{(c)} = 0,
\end{equation}
which implies that
\begin{equation}
    \boldsymbol{d}_{j + \frac{1}{2}}^{(p)} = \boldsymbol{d}_{j - \frac{1}{2}}^{(p)} = 0.
\end{equation}

Using \eqref{average flux at stationary contact}, the averaged fluxes appearing in \eqref{j + 1/2 flux} and \eqref{j - 1/2 flux} both reduce to
\begin{equation}
    \hat{\textbf{f}}_{j + \frac{1}{2}} = \hat{\textbf{f}}_{j - \frac{1}{2}} = \begin{bmatrix}
        0 & p_0 & 0
    \end{bmatrix}^T.
\end{equation}
Therefore, the finite volume update reduces to \eqref{discrete equilibrium} and the modified PVU scheme preserves the stationary contact discontinuity as a discrete equilibrium.


\section{A compactly supported distribution function in two dimensions}
\label{sec 6}

\subsection{Moments of the Maxwellian distribution function in 2D}
\label{subsec 6.1}

In this subsection, we study the moments of the two-dimensional Maxwell velocity distribution function. We first consider the moments that recover the conserved variable vector of the two-dimensional Euler equations of gas dynamics. The zeroth moment of the Maxwellian distribution function gives the first component of the conserved variable vector, \textit{i.e.},
\begin{equation*}
    U_1 = \int_{0}^{\infty} \, dI \int_{-\infty}^{\infty} \, dc_2 \int_{-\infty}^{\infty} \, dc_1 \, {f}^{eq}
\end{equation*}
Integrating the above equation with respect to I leads to
\begin{equation} \label{2D_zeroth_F_moment}
    \int_{-\infty}^{\infty} \, dc_2 \int_{-\infty}^{\infty} \, dc_1 \, \Bar{f}^{eq} = \rho,
\end{equation}
where $ \Bar{f}^{eq} $ denotes the Maxwellian distribution function integrated with respect to $ I $, which reads
\begin{equation}
    \Bar{f}^{eq} = \rho \frac{\beta}{\pi} \exp\biggl( -\beta \left( c_1^2 + c_2^2 \right) \biggl).
\end{equation}
The first moment of the Maxwellian distribution function with respect to the x component of the bulk fluid velocity $ u_1 $ and the peculiar velocity $ c_1 $ recovers the second component of the conserved variable vector, \textit{i.e.},
\begin{equation*}
    U_2 = \int_0^{\infty} \, dI \int_{-\infty}^{\infty} \, dc_2 \int_{-\infty}^{\infty} \, dc_1 \, \left( u_1 + c_1 \right) f^{eq} = \rho u_1.
\end{equation*}
Integrating the above equation with respect to $ I $ reads
\begin{equation} \label{2D_first_1_F_moment}
    \int_{-\infty}^{\infty} \, dc_2 \int_{-\infty}^{\infty} \, dc_1 \, c_1 \Bar{f}^{eq} = 0.
\end{equation}
Similarly, the first moment of the Maxwellian distribution function corresponding to the y component of the bulk fluid velocity $ u_2 $ and the peculiar velocity $ c_2 $ gives the third component of the conserved variable vector, which reads
\begin{equation*}
    U_3 = \int_0^{\infty} \, dI \int_{-\infty}^{\infty} \, dc_2 \int_{-\infty}^{\infty} \, dc_1 \, \left( u_2 +c_2 \right) f^{eq} = \rho u_2.
\end{equation*}
Integration of the above equation with respect to $ I $ gives
\begin{equation} \label{2D_first_2_F_moment}
    \int_{-\infty}^{\infty} \, dc_1 \int_{-\infty}^{\infty} \, dc_2 \, c_2 \Bar{f}^{eq} = 0.
\end{equation}
The moment of the Maxwellian velocity distribution function corresponding to the fourth component of the conserved variable vector reads
\begin{equation*}
    U_4 =\int_0^{\infty} \, dI \int_{-\infty}^{\infty} \, dc_2 \int_{-\infty}^{\infty} \, dc_1 \Bigl(I + \dfrac{u_1^2}{2} + \dfrac{c_1^2}{2} + u_1 c_1 + \dfrac{u_2^2}{2} + \dfrac{c_2^2}{2} + u_2 c_2 \Bigl) f^{eq}.
\end{equation*}
Integrating the above equation with respect to $ I $ gives
\begin{equation*}
    \int_{-\infty}^{\infty} \, dc_2 \int_{-\infty}^{\infty} \, dc_1 \, \bigl( c_1^2 + c_2^2 \bigl) \Bar{f}^{eq} = 2p.
\end{equation*}
Further simplification of the above equation leads to the following moment relation:
\begin{equation} \label{eq:2D_first_3_F_moment}
    \int_{-\infty}^{\infty} \, dc_2 \int_{-\infty}^{\infty} \, dc_1 \, c_1^2 \Bar{f}^{eq} = p.
\end{equation}
Next, we consider the moments of the Maxwellian distribution function corresponding to the flux vectors of the two-dimensional Euler equations of gas dynamics. We note that the moments integrated with respect to $ I $ corresponding to the first component of the flux vectors in the x and y directions give the same expressions as mentioned in equations (\ref{2D_first_1_F_moment}) and (\ref{2D_first_2_F_moment}). The moment concerning the second component of the flux vector in the x direction reproduces the same moment relation as obtained in the equation (\ref{eq:2D_first_3_F_moment}). The next moment corresponds to the third component of the flux vector in the x direction and the second component of the flux vector in the y direction, which reads
\begin{equation*}
    G_{13} = G_{22} = \int_0^{\infty} \, dI \int_{-\infty}^{\infty} \, dc_2 \int_{-\infty}^{\infty} \, dc_1 \, \left( u_1 + c_1 \right) \left( u_2 + c_2 \right) f^{eq}.
\end{equation*}
Integrating the above equation with respect to $ I $ leads to
\begin{equation}
    \int_{-\infty}^{\infty} \, dc_2 \int_{-\infty}^{\infty} \, dc_1 \, \left( u_1 + c_1 \right) \left( u_2 + c_2 \right) \Bar{f}^{eq} = \rho u_1 u_2.
\end{equation}
Further simplification of the above equation leads to
\begin{equation} \label{2D_second_13_F_moment}
    \int_{-\infty}^{\infty} \, dc_2 \int_{-\infty}^{\infty} \, dc_1 \, c_2 \, c_1 \Bar{f}^{eq} = 0.
\end{equation}
The moment corresponding to the third component of the flux vector in the y direction reads
\begin{equation*}
    G_{23} = \int_0^{\infty} \, dI \int_{-\infty}^{\infty} \, dc_2 \int_{-\infty}^{\infty} \, dc_1 \, \left( u_2 + c_2 \right)^2 f^{eq}.
\end{equation*}
Integrating the above equation with respect to $ I $ gives the following moment relation: 
\begin{equation} \label{2D_second_23_F_moment}
    \int_{-\infty}^{\infty} \, dc_1 \int_{-\infty}^{\infty} \, dc_2 \, c_2^2 \, \Bar{f}^{eq} = p.
\end{equation}
The moment corresponding to the fourth component of the flux vector in the x direction is written as
\begin{equation*}
     G_{14} = \int_0^{\infty} \, dI \int_{-\infty}^{\infty} \, dc_2 \int_{-\infty}^{\infty} \, dc_1 \, \Bigl(I + \dfrac{\left( u_1 + c_1 \right)^2}{2} + \dfrac{\left( u_2 + c_2 \right)^2}{2} \Bigl) \left( u_1 + c_1 \right) f^{eq}.
\end{equation*}
Integrating the above equation with respect to $ I $ and rearranging the terms lead to
\begin{equation} \label{2D_third_1_F_moment}
    \int_{-\infty}^{\infty} \, dc_2 \int_{-\infty}^{\infty} \, dc_1 \, c_1^3 \, \Bar{f}^{eq} = 0.
\end{equation}
The moment concerning the fourth component of the flux vector in the y direction reads
\begin{equation*}
    G_{24} = \int_0^{\infty} \, dI \int_{-\infty}^{\infty} \, dc_2 \int_{-\infty}^{\infty} \, dc_1 \, \Bigl(I + \dfrac{\left( u_1 + c_1 \right)^2}{2} + \dfrac{\left( u_2 + c_2 \right)^2}{2} \Bigl) \left( u_2 + c_2 \right) f^{eq}.
\end{equation*}
Integrating the above equation with respect to $ I $ and rearranging the terms give
\begin{equation} \label{2D_third_2_F_moment}
    \int_{-\infty}^{\infty} \, dc_1 \int_{-\infty}^{\infty} \, dc_2 \, c_2^3 \, \Bar{f}^{eq} = 0.
\end{equation}
The study of all the twelve moments corresponding to the conserved variable and flux vectors of the two-dimensional Euler equations of the gas dynamics reveals that there are only eight independent moment relations. Any distribution function $ \Bar{f}^{eq} $ that satisfies these eight integrals will recover the Euler equations in two dimensions for $ \gamma = 5/3 $. We list only the three nonzero independent moment relations together for future reference.

\begin{subequations} \label{2D_reduced_F_moments}
    \begin{align} 
        \int_{-\infty}^{\infty} \, dc_2 \int_{-\infty}^{\infty} \, dc_1 \, \Bar{f}^{eq} &= \rho, \label{2D_reduced_zeroth_F_moment} \\
        \int_{-\infty}^{\infty} \, dc_2 \int_{-\infty}^{\infty} \, dc_1 \, c_1^2 \, \Bar{f}^{eq} &= p, \label{2D_reduced_second_v1F_moment} \\
        \int_{-\infty}^{\infty} \, dc_1 \int_{-\infty}^{\infty} \, dc_2 \, c_2^2 \, \Bar{f}^{eq} &= p. \label{2D_reduced_second_v2F_moment}
    \end{align}
\end{subequations}

\subsection{A compactly supported hat function in two dimensions}
\label{subsec 6.2}

\begin{figure}[h]
   \centering
   \includegraphics[scale=0.6]{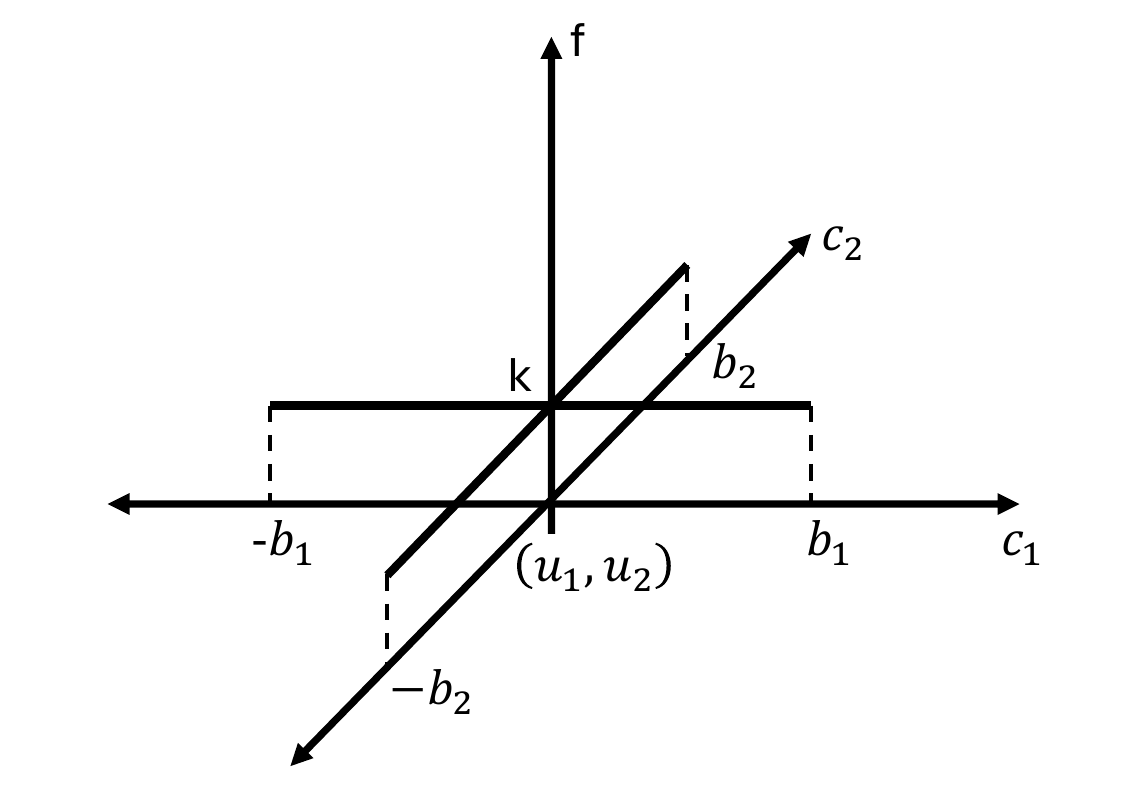}
   \caption{A two-dimensional hat function}
   \label{fig:2D_F_approximation}
\end{figure}

In this subsection, we construct the compactly supported distribution function in two dimensions that satisfies the three nonzero independent moment relations listed in Equation (\ref{2D_reduced_F_moments}). We replace the two-dimensional Maxwellian distribution function with a product of hat functions in $ c_1 $ and $ c_2 $, as shown in figure \ref{fig:2D_F_approximation}.
The new distribution function then reads
\begin{equation} \label{2D_f^eq}
    \Bar{f}^{eq} = \rho f^{(x)}(c_1) f^{(y)}(c_2),
\end{equation}
where $ f^{(x)}(c_1) $ and $ f^{(y)}(c_2) $ are defined as
\begin{equation} \label{New_2D_xDistribution_function}
    f^{(x)}(c_1)  = 
    \begin{cases}
        0 & \quad \text{if} \quad c_1 < -b_1 \\[7 pt]
        k & \quad \text{if} \quad -b_1 \leq c_1 \leq b_1 \\[7 pt]
        0 & \quad \text{if} \quad c_1 > b_1 \\
    \end{cases},
\end{equation}

\begin{equation} \label{New_2D_yDistribution_function}
    f^{(y)}(c_2)  = 
    \begin{cases}
        0 & \quad \text{if} \quad c_2 < -b_2 \\[7 pt]
        k & \quad \text{if} \quad -b_2 \leq c_2 \leq b_2 \\[7 pt]
        0 & \quad \text{if} \quad c_2 > b_2 \\
    \end{cases}.
\end{equation}
The maximum value $ k $ and the half-widths $ b_1 $ and $ b_2 $ of the new distribution function are determined by enforcing that it satisfies three nonzero independent moment relations. Substituting the new distribution function into Equation (\ref{2D_reduced_zeroth_F_moment}), we obtain
\begin{equation} \label{2D_Fb1b2_condition_1}
    4\, k^2 \, b_1 \, b_2 = 1.
\end{equation}
We then substitute the new ansatz into Equation (\ref{2D_reduced_second_v1F_moment}) to obtain
\begin{equation} \label{2D_b1_condition_2}
    \dfrac{4}{3} \, k^2 \, b_1^3 \, b_2 = R \, T.
\end{equation} 
Finally, substituting the new distribution function into equation (\ref{2D_reduced_second_v2F_moment}) to get
\begin{equation} \label{2D_b2_condition_3}
    \dfrac{4}{3} \, k^2 \, b_1 \, b_2^3 = R \, T.
\end{equation}
Using Equation (\ref{2D_Fb1b2_condition_1}) in Equations (\ref{2D_b1_condition_2}) and (\ref{2D_b2_condition_3}) leads to
\begin{equation} \label{b1_b2}
    b_1 = \sqrt{3 R T}, \; b_2 = \sqrt{3 R T}.
\end{equation}
Using Equation (\ref{b1_b2}) in (\ref{2D_Fb1b2_condition_1}) to obtain
\begin{equation}
    k = \dfrac{1}{2 \sqrt{3 R T}}.
\end{equation}
\textbf{Remark:} The compactly supported distribution function shown in (\ref{2D_f^eq}), (\ref{New_2D_xDistribution_function}), (\ref{New_2D_yDistribution_function}) will satisfy the relations (\ref{2D_first_1_F_moment}), (\ref{2D_first_2_F_moment}), (\ref{2D_second_13_F_moment}), (\ref{2D_third_1_F_moment}) and (\ref{2D_third_2_F_moment}) identically.\\
Thus, the compactly supported hat function in two dimensions is defined as
\begin{equation} \label{2D_new_distribution_function}
    \Bar{f}^{eq} = \rho \, f^{(x)}(c_1) \, f^{(y)}(c_2),
\end{equation}
where 
\begin{equation} \label{2D_xlinearDistribution_function}
    f^{(x)}(c_1) = 
    \begin{cases}
        \; \quad 0 & \quad \text{if} \quad c_1 < -\sqrt{3 R T} \\[7 pt]
        \dfrac{1}{2 \sqrt{3 R T}} & \quad \text{if} \quad -\sqrt{3 R T} \leq c_1 \leq \sqrt{3 R T} \\[7 pt]
        \; \quad 0 & \quad \text{if} \quad c_1 > \sqrt{3 R T} \\
    \end{cases},
\end{equation}
\begin{equation} \label{2D_ylinearDistribution_function}
    f^{(y)}(c_2) = 
        \begin{cases}
        \; \quad 0 & \quad \text{if} \quad c_2 < -\sqrt{3 R T} \\[7 pt]
        \dfrac{1}{2 \sqrt{3 R T}} & \quad \text{if} \quad -\sqrt{3 R T} \leq c_2 \leq \sqrt{3 R T} \\[7 pt]
        \; \quad 0 & \quad \text{if} \quad c_2 > \sqrt{3 R T} \\
    \end{cases}.
\end{equation}

\subsection{Equivalent form of the Euler equations of gas dynamics in 2D}
\label{subsec 6.3}

The moments of the compactly supported equilibrium distribution function defined in (\ref{2D_new_distribution_function}), (\ref{2D_xlinearDistribution_function}), (\ref{2D_ylinearDistribution_function}) recover the conserved variable and the flux vectors of the Euler equations for monatomic gases, \textit{i.e.},
\begin{subequations}
    \begin{align}
        \int_{\boldsymbol{c}} \, D \boldsymbol{c} \, \Tilde{\boldsymbol{\psi}} \bar{f}^{eq} &= \begin{bmatrix}
            \rho & \rho u_1 & \rho u_2 & \dfrac{1}{2} \rho | \boldsymbol{u} |^2 + p
        \end{bmatrix}^T, \\
        \int_{\boldsymbol{c}} \, D \boldsymbol{c} \, \Tilde{\boldsymbol{\psi}} \left( u_1 + c_1 \right) \bar{f}^{eq} &= \begin{bmatrix}
            \rho u_1 & p + \rho u_1^2 & \rho u_1 u_2 & \left( \dfrac{1}{2} \rho | \boldsymbol{u} |^2 + 2 p \right) u_1
        \end{bmatrix}^T, \\
        \int_{\boldsymbol{c}} D \boldsymbol{c} \, \Tilde{\boldsymbol{\psi}} \left( u_2 + c_2 \right) \bar{f}^{eq} &= \begin{bmatrix}
            \rho u_2 & \rho u_1 u_2 & p + \rho u_2^2 & \left( \dfrac{1}{2} \rho | \boldsymbol{u} |^2 + 2 p \right) u_2
        \end{bmatrix}^T,
    \end{align}
\end{subequations}
where $ \Tilde{\boldsymbol{\psi}} $ is defined as
\begin{equation}
    \Tilde{\boldsymbol{\psi}} = \begin{bmatrix}
        1 & \left( u_1 + c_1 \right) & \left( u_2 + c_2 \right) & \dfrac{1}{2} | \boldsymbol{u} + \boldsymbol{c} |^2
    \end{bmatrix}^T.
\end{equation}
To take into account the internal degrees of freedom of molecules, Perthame proposed an approach in \cite{Perthame1999}, introducing two compactly supported distribution functions, $ \bar{f}^{eq} $ and $ \bar{g}^{eq} $, such that
\begin{equation}
    \int_{\boldsymbol{c}} \, D \boldsymbol{c} \, \dfrac{| \boldsymbol{u} + \boldsymbol{c} |^2}{2} \bar{f}^{eq} + \int_{\boldsymbol{c}} \, D \boldsymbol{c} \, \bar{g}^{eq} = \dfrac{1}{2} \rho \boldsymbol{u}^2 + \dfrac{p}{\left( \gamma - 1 \right)}.
\end{equation}
A simple choice for $ \bar{g}^{eq} $ is
\begin{equation} \label{2D_second_distribution_function}
    \bar{g}^{eq} = I_0 \bar{f}^{eq},
\end{equation}
where $ I_0 $ is defined in (\ref{definition_of_I_0}) with $ D_f = 2 $. With this choice, the conserved variable and the flux vectors corresponding to the Euler equations for polyatomic gases can be expressed as
\begin{subequations}
    \begin{align}
        \boldsymbol{U} &= \int_{\boldsymbol{c}} \, D \boldsymbol{c} \, \Tilde{\boldsymbol{\psi}}_1 \, \bar{f}^{eq} + \int_{\boldsymbol{c}} \, D \boldsymbol{c} \, \Tilde{\boldsymbol{\psi}}_2 \, \bar{g}^{eq}, \label{2D_definition_of_conserved_variable} \\
        \boldsymbol{G}_1 &= \int_{\boldsymbol{c}} \, D \boldsymbol{c} \, \Tilde{\boldsymbol{\psi}}_1 \, \left( u_1 + c_1 \right) \, \bar{f}^{eq} + \int_{\boldsymbol{c}} \, D \boldsymbol{c} \, \Tilde{\boldsymbol{\psi}}_2 \, \left( u_1 + c_1 \right) \, \bar{g}^{eq}, \\
        \boldsymbol{G}_2 &= \int_{\boldsymbol{c}} \, D \boldsymbol{c} \, \Tilde{\boldsymbol{\psi}}_1 \, \left( u_2 + c_2 \right) \, \bar{f}^{eq} + \int_{\boldsymbol{c}} \, D \boldsymbol{c} \, \Tilde{\boldsymbol{\psi}}_2 \, \left( u_2 + c_2 \right) \, \bar{g}^{eq},
    \end{align}
\end{subequations}
where $ \Tilde{\boldsymbol{\psi}}_1 $ and $ \Tilde{\boldsymbol{\psi}}_2 $ are defined as
\begin{equation}
    \Tilde{\boldsymbol{\psi}}_1 = \begin{bmatrix}
        1 & \left( u_1 + c_1 \right) & \left( u_2 + c_2 \right) & \dfrac{1}{2} |\boldsymbol{u} + \boldsymbol{c}|^2
    \end{bmatrix}^T, \; \Tilde{\boldsymbol{\psi}}_2 = \begin{bmatrix}
        0 & 0 & 0 & 1
    \end{bmatrix}^T.
\end{equation}
The flow is now described by two velocity distribution functions, $ f(\boldsymbol{x}, \boldsymbol{u}, \boldsymbol{c}, t) $ and $ g(\boldsymbol{x}, \boldsymbol{u}, \boldsymbol{c}, t) $, at spatial location $ \boldsymbol{x} $ and time $ t $, with macroscopic velocity $ \boldsymbol{u} $ and peculiar velocity $ \boldsymbol{c} $. These functions evolve according to the following Boltzmann-like equations:
\begin{subequations} \label{2D_Boltzmann-like_equations}
    \begin{align}
        \dfrac{\partial f}{\partial t} + \nabla \cdot \left( \boldsymbol{u} f \right) + \nabla \cdot \left( \boldsymbol{c} f \right) = Q_1, \\
        \dfrac{\partial g}{\partial t} + \nabla \cdot \left( \boldsymbol{u} g \right) + \nabla \cdot \left( \boldsymbol{c} g \right) = Q_2,
    \end{align}
\end{subequations}
where $ Q_1 $ and $ Q_2 $ are unknown functions satisfying for all $ t $, $ x $

\begin{subequations}
    \begin{align}
        \int_{\boldsymbol{c}} \, D \boldsymbol{c} \, Q_1 = 0, \\
        \int_{\boldsymbol{c}} \, D \boldsymbol{c} \, \left( u_1 + c_1 \right) \, Q_1 = 0, \\
        \int_{\boldsymbol{c}} \, D \boldsymbol{c} \, \left( u_2 + c_2 \right) \, Q_1 = 0, \\
        \int_{\boldsymbol{c}} \, D \boldsymbol{c} \, \left( \dfrac{| \boldsymbol{u} + \boldsymbol{c} |^2}{2} \, Q_1 + Q_2 \right) = 0.
    \end{align}
\end{subequations}
Thus, in the limits $ f \rightarrow \bar{f}^{eq} $ and $ g \rightarrow \bar{g}^{eq} $, the two-dimensional Euler equations of gas dynamics can be recovered as a linear combination of the two Boltzmann-like equations (\ref{2D_Boltzmann-like_equations}), \textit{i.e.},
\begin{multline}
    \int_{\boldsymbol{c}} \, D \boldsymbol{c} \, \Tilde{\boldsymbol{\psi}}_1 \left( \dfrac{\partial f}{\partial t} + \nabla \cdot \left( \boldsymbol{u} f \right) + \nabla \cdot \left( \boldsymbol{c} f \right), \; f = \bar{f}^{eq} \right) \\
    + \int_{\boldsymbol{c}} \, D \boldsymbol{c} \, \Tilde{\boldsymbol{\psi}}_2 \left( \dfrac{\partial g}{\partial t} + \nabla \cdot \left( \boldsymbol{u} g \right) + \nabla \cdot \left( \boldsymbol{c} g \right), \; g = \bar{g}^{eq} \right) = 0.
\end{multline}


\section{A contact discontinuity capturing PVU Boltzmann scheme based on hat function in 2D}
\label{sec 7}

\subsection{Numerical formulation}
\label{subsec 7.1}
Let $ T_h $ be a tessellation of the $ \Re^2 $ space. Let $ m(K) $ be a measure in $ \Re^2 $ of the cell $ K $, $ e $ be a side of $ K $, $ K_e $ be a neighboring cell of $ K $ along the edge $ e $, $ \hat{\textbf{n}}_e $ be an outward unit normal vector to the edge $ e $ making an angle $ \theta_e $ with the $ x- $axis, $ \hat{\textbf{t}}_e $ be a unit tangent vector to the edge, $ m(e) $ be a measure in $ \Re $ of the edge $ e $. Initially, the fluid velocity vector and the thermodynamic state variables are assumed to be piecewise constant in each finite volume cell. Consequently, the equilibrium distribution functions defined in (\ref{2D_new_distribution_function}) and (\ref{2D_second_distribution_function}) are constructed locally within each cell. The distribution functions in each control volume evolve in time by solving Boltzmann-like equations (\ref{2D_Boltzmann-like_equations}) using a two-stage procedure. \\

\textbf{Convection step:} We start with the two-dimensional collision-less peculiar velocity based Boltzmann equations, which reads

\begin{subequations} \label{2D_collisionless_PV_Boltzmann_equation}
    \begin{align}
        \frac{\partial f}{\partial t} + \nabla \cdot \left( \boldsymbol{u} f \right) + \nabla \cdot \left( \boldsymbol{c} f \right) = 0, \\
        \frac{\partial g}{\partial t} + \nabla \cdot \left( \boldsymbol{u} g \right) + \nabla \cdot \left( \boldsymbol{c} g \right) = 0.
    \end{align}
\end{subequations}
Integrating (\ref{2D_collisionless_PV_Boltzmann_equation}) over the cell area and applying Gauss's divergence theorem yields
\begin{subequations}
    \begin{align}
        \int_{K} \, dA \, \dfrac{\partial f}{\partial t} = -\oint \, dl \, \hat{f}, \\
        \int_{K} \, dA \, \dfrac{\partial g}{\partial t} = -\oint \, dl \, \hat{g},
    \end{align}
\end{subequations}
where
\begin{equation}
    \hat{h} = \left( \boldsymbol{u} \cdot \hat{\boldsymbol{n}} + \boldsymbol{c} \cdot \hat{\boldsymbol{n}} \right) h, 
\end{equation}
for $ h = $ $ f $ or $ g $. The total flux through the cell boundary can be expressed as the sum of the one-dimensional fluxes evaluated at the midpoints of each line segment, \textit{i.e.},
\begin{subequations}
    \begin{align}
        \int_{K} \, dA \, \dfrac{\partial f}{\partial t} = -\sum_e \hat{f}_e \, m(e), \\
        \int_{K} \, dA \, \dfrac{\partial g}{\partial t} = -\sum_e \hat{g}_e \, m(e),
    \end{align}
\end{subequations}
with
\begin{equation}
    \hat{h}_e = \left( \boldsymbol{u}_e \cdot \hat{\boldsymbol{n}}_e + \boldsymbol{c}_e \cdot \hat{\boldsymbol{n}}_e \right) h_e.
\end{equation}
The forward-in-time discretization of the temporal derivative leads to
\begin{subequations} \label{discretized_2D_PV_Boltzmann_equation}
    \begin{align}
        \int_{K} \, dA \, f_K^{(n + 1)-} = \int_{K} \, dA \, f_K^{(n)} - \Delta t \, \sum_e \hat{f}^{(n)}_e \, m(e), \\
        \int_{K} \, dA \, g_K^{(n + 1)-} = \int_{K} \, dA \, g_K^{(n)} - \Delta t \, \sum_e \hat{g}^{(n)}_e \, m(e).
    \end{align}
\end{subequations}

\textbf{Collision step:} The collision step restores equilibrium in each finite volume cell, \textit{i.e.},
\begin{subequations} \label{essence_of_collision_step_2D}
    \begin{align}
        \int_K \, dA \, \int_{\boldsymbol{c}} \, D \boldsymbol{c} \, \Tilde{\boldsymbol{\psi}}_1 \, \bar{f}^{eq, (n + 1)}_K &= \int_K \, dA \, \int_{\boldsymbol{c}} \, D \boldsymbol{c} \, \Tilde{\boldsymbol{\psi}}_1 \, f^{(n + 1)-}_K, \\
        \int_K \, dA \, \int_{\boldsymbol{c}} \, D \boldsymbol{c} \, \Tilde{\boldsymbol{\psi}}_2 \, \bar{g}^{eq, (n + 1)}_K &= \int_K \, dA \, \int_{\boldsymbol{c}} \, D \boldsymbol{c} \, \Tilde{\boldsymbol{\psi}}_2 \, g^{(n + 1)-}_K.
    \end{align}
\end{subequations}
Substituting Equations (\ref{discretized_2D_PV_Boltzmann_equation}) and (\ref{essence_of_collision_step_2D}) into the definition (\ref{2D_definition_of_conserved_variable}) yields a first-order finite volume scheme for the two-dimensional Euler equations of gas dynamics, which reads
\begin{equation} \label{2D_finite_volume_update_formula}
    \bar{\boldsymbol{U}}_K^{n + 1} = \bar{\boldsymbol{U}}_K^{n} - \sigma \sum_e \hat{\textbf{f}}_e,
\end{equation}
where $ \bar{\boldsymbol{U}} $ denotes the cell-averaged value of the conserved variable vector in each cell, $ \hat{\textbf{f}}_e $ denotes the numerical flux in the direction normal to the edge $ e $, which is represented by the sum of the numerical advective and pressure fluxes, \textit{i.e.},
\begin{equation} \label{2D_definition_of_total_flux}
    \hat{\textbf{f}}_e = \hat{\textbf{f}}^{(a)}_e + \hat{\textbf{f}}^{(p)}_e.
\end{equation}
We choose
\begin{multline} \label{2D_definition_advective_numerical_flux_PVU_plus}
    \hat{\textbf{f}}_{e}^{(a)} = \dfrac{1}{2} \int_{c_t} \, dc_t \int_{c_n} \, dc_n \, \Tilde{\boldsymbol{\psi}}_{1, e} \left( (u_n \bar{f}^{eq})_K + (u_n \bar{f}^{eq})_{K_e} - \Delta_e \left( \phi^{(u)}_e |u_n| \bar{f}^{eq} \right) \right) \\
    + \dfrac{1}{2} \int_{c_t} \, dc_t \int_{c_n} \, dc_n \, \Tilde{\boldsymbol{\psi}}_{2, e} \left( (u_n \bar{g}^{eq})_K + (u_n \bar{g}^{eq})_{K_e} - \Delta_{e} \left( \phi^{(u)}_{e} |u_n| \bar{g}^{eq} \right) \right),
\end{multline}
\begin{multline} \label{2D_definition_pressure_numerical_flux_PVU_plus}
    \hat{\textbf{f}}_{e}^{(p)} = \dfrac{1}{2} \int_{c_t} \, dc_t \int_{c_n} \, dc_n \, \Tilde{\boldsymbol{\psi}}_{1, e} \left( (c_n \bar{f}^{eq})_K + (c_n \bar{f}^{eq})_{K_e} - \Delta_e \left( \phi^{(c)}_e |c_n| \bar{f}^{eq} \right) \right) \\
    + \dfrac{1}{2} \int_{c_t} \, dc_t \int_{c_n} \, dc_n \, \Tilde{\boldsymbol{\psi}}_{2, e} \left( (c_n \bar{g}^{eq})_K + (c_n \bar{g}^{eq})_{K_e} - \Delta_{e} \left( \phi^{(c)}_{e} |c_n| \bar{g}^{eq} \right) \right),
\end{multline}
where
\begin{equation} \label{definition_of_delta_e}
    \Delta_{e} (\cdot) = (\cdot)_{K_e} - (\cdot)_{K},
\end{equation}
and
\begin{subequations} \label{normal_bulk_fluid_and_peculiar_velocity}
    \begin{align}
        u_n &= \textbf{u} \cdot \hat{\textbf{n}}_e, \, u_t = \textbf{u} \cdot \hat{\textbf{t}}_e, \\
        c_n &= \textbf{c} \cdot \hat{\textbf{n}}_e, \, c_t = \textbf{c} \cdot \hat{\textbf{t}}_e.
    \end{align}
\end{subequations}
$ \Tilde{\boldsymbol{\psi}}_{1, e} $ and $ \Tilde{\boldsymbol{\psi}}_{2, e} $ denote the moment vectors expressed in the local coordinate system aligned with the normal and tangential directions, \textit{i.e.},
\begin{equation} \label{2D_defintion_of_moment_vector_1}
    \Tilde{\boldsymbol{\psi}}_{1, e} = \begin{bmatrix}
        1 & v_n \cos \theta_e - v_t \sin \theta_e & v_n \sin \theta_e + v_t \cos \theta_e & \dfrac{1}{2} (v_n^2 + v_t^2)
    \end{bmatrix}^T,
\end{equation}
where
\begin{subequations}
    \begin{align}
        v_n &= u_n + c_n, \\
        v_t &= u_t + c_t.
    \end{align}
\end{subequations}
and
\begin{equation} \label{2D_defintion_of_moment_vector_2}
    \Tilde{\boldsymbol{\psi}}_{2, e} = \begin{bmatrix}
        0 & 0 & 0 & 1
    \end{bmatrix}^T.
\end{equation}
The diffusion control parameters $ \phi^{(u)}, \phi^{(c)} $ obey the dynamical consistency condition, \textit{i.e.},
\begin{equation}
    \textbf{U}_K = \textbf{U}_{K_e} \implies \left( \phi^{(u), (c)}_e \right)_K = \left( \phi^{(u), (c)}_e \right)_{K_e}.
\end{equation}
Furthermore, the diffusion control parameters $ \phi^{(u)}_e $ and $ \phi^{(c)}_e $ satisfy the constraints summarized in Equation (\ref{constraints_on_phi's}).

\subsection{Positivity analysis in two dimensions}
\label{subsec 7.2}
In this subsection, we prove the positivity-preserving property of the 2D modified PVU Boltzmann scheme while keeping the $ \phi $'s general.

\vspace{1\baselineskip}
\noindent
\textbf{Theorem 8.} Assume $ \Bar{f}^{eq}_K = 0 $ for $ |(c_n)_K| \geq |(c_n)_K|_{max} $. Furthermore, assume that the diffusion control parameters $ \phi^{(u)}_e, \phi^{(c)}_e $ obey the following constraints:
\begin{subequations}
    \begin{align}
        \dfrac{1}{2} \sum_{e} |(u_n)_K| \, \phi^{(u)}_e \, m(e) &\leq m_1 \displaystyle \max_{\forall K} \{ |(u_n)_K| \} \sum_{e} m(e), \label{additional_constraint_on_phi1_2D} \\ 
        \dfrac{1}{2} \sum_{e} |(c_n)_K|_{max} \, \phi^{(c)}_e \, m(e) &\leq m_2 \displaystyle \max_{\forall K} \{ |(c_n)_K|_{max} \} \sum_{e} m(e), \label{additional_constraint_on_phi2_2D}
    \end{align}
\end{subequations}
where $ m_1, m_2 \in \mathbb{N} $. Then, under the CFL-like condition,
\begin{equation} \label{2D_CFL_like_condition}
    \frac{\Delta t}{m(K)} \max_{K} \{ m_1 \, |(u_n)_K| + m_2 \, |(c_n)_K|_{max} \} \sum_{e} m(e) \leq 1,
\end{equation}
the two-dimensional modified PVU Boltzmann scheme is positive, \textit{i.e.}, $ \rho > 0 $ and $ p > 0 $.

\vspace{1\baselineskip}
\noindent
Proof: We have, from (\ref{2D_finite_volume_update_formula}), an update formula for the density, which reads
\begin{equation}
    \rho_K^{(n + 1)} = \left\langle \Bar{f}_K^{eq, (n)} - \frac{\Delta t}{m(K)} \sum_{e} \hat{\textbf{f}}_e \, m(e) \right\rangle.
\end{equation}
We have, from (\ref{2D_definition_advective_numerical_flux_PVU_plus}) and (\ref{2D_definition_pressure_numerical_flux_PVU_plus}),
\begin{multline} \label{2D_update_formula_for_density}
    \rho_K^{(n + 1)} = \left\langle \Bar{f}_K^{eq, (n)} \left( 1 - \frac{\Delta t}{2 m(K)} \sum_{e} \left( \phi^{(u)}_e |(u_n)_K| + \phi^{(c)}_e |(c_n)_K| \right) m(e) \right) \right\rangle \\
    - \left\langle \frac{\Delta t}{2 m(K)} \sum_{e} h_{K_e}^- m(e) \right\rangle,
\end{multline}
where
\begin{equation}
    h_{K_e}^- = \left( u_n - |u_n| \phi^{(u)}_e \right)_{K_e} \Bar{f}^{eq}_{K_e} + \left( c_n - |c_n| \phi^{(c)}_e \right)_{K_e} \Bar{f}^{eq}_{K_e}.
\end{equation}
Equation (\ref{2D_update_formula_for_density}) is possible in view of the following relations:
\begin{equation}
    \textbf{u}_K \cdot \sum_{e} \hat{\textbf{n}}_e \, m(e) = 0, \; \textbf{c}_K \cdot \sum_{e} \hat{\textbf{n}}_e \, m(e) = 0.
\end{equation}
We then note from Proposition 2. and (\ref{positivity_of_u_tilda_plus}), (\ref{positivity_of_u_tilda_minus}) that
\begin{equation}
    -\left\langle \frac{\Delta t}{2 m(K)} \sum_{e} h_{K_e}^- m(e) \right\rangle \geq 0.
\end{equation}
We only have to consider such peculiar velocities for which $ \Bar{f}^{eq} > 0 $, \textit{i.e.}, 
\begin{equation}
    |(c_n)_K| \leq |(c_n)_K|_{max} \leq \max_{\forall K} \{ |(c_n)_K|_{max} \}.
\end{equation}
We have, from (\ref{additional_constraint_on_phi2_2D}),
\begin{equation}
    \sum_{e} |(c_n)_K| \, \phi^{(c)}_e \, m(e) \leq \, 2 \, m_2 \, \max_{\forall K} \, \{ |(c_n)_K|_{max} \} \, \sum_{e} m(e).
\end{equation}
We note, from (\ref{additional_constraint_on_phi1_2D}),
\begin{equation}
    \sum_{e} |(u_n)_K| \, \phi^{(u)}_e \, m(e) \leq 2 \, m_1 \max_{\forall K} \{ |(u_n)_K| \} \, \sum_{e} m(e).
\end{equation}
We add the above two inequalities to obtain
\begin{multline} \label{final_step_in_positivity_proof_of_density}
    \sum_{e} \left( \phi^{(u)}_e |(u_n)_K| + \phi^{(c)}_e |(c_n)_K| \right) m(e) \\
    \leq 2 \max_{K} \{ m_1 \, |(u_n)_K| + m_2 \, |(c_n)_K|_{max} \} \sum_{e} m(e).
\end{multline}
Then the CFL-like condition (\ref{2D_CFL_like_condition}) dictates that the first term in (\ref{2D_update_formula_for_density}) is positive and hence
\begin{equation}
    \rho_K^{(n+1)} \geq 0.
\end{equation}
We now show that the 2D modified PVU Boltzmann scheme ensures the positivity of pressure.
We begin with an update formula for pressure:
\begin{equation}
    \begin{aligned}
        p_K^{(n + 1)} = &\left\langle c_K^2 \bar{f}_K^{eq, (n)} \left( 1 - \frac{\Delta t}{2 m(K)} \sum_{e} \left( \phi^{(u)}_e |(u_n)_K| + \phi^{(c)}_e |(c_n)_K| \right) m(e) \right) \right\rangle \\
        &- \left\langle \frac{\Delta t}{2 m(K)} \sum_{e} c_K^2 h_{K_e}^- m(e) \right\rangle,
    \end{aligned}
\end{equation}
We note from Proposition 3. and (\ref{positivity_of_u_tilda_plus}), (\ref{positivity_of_u_tilda_minus}) that
\begin{equation}
    - \left\langle \frac{\Delta t}{2 m(K)} \sum_{e} c_K^2 h_{K_e}^- m(e) \right\rangle \geq 0.
\end{equation}
We further note, from (\ref{final_step_in_positivity_proof_of_density}), that
\begin{equation}
    \left\langle c_K^2 \bar{f}_K^{eq, (n)} \left( 1 - \frac{\Delta t}{2 m(K)} \sum_{e} \left( \phi^{(u)}_e |(u_n)_K| + \phi^{(c)}_e |(c_n)_K| \right) m(e) \right) \right\rangle \geq 0.
\end{equation}
Hence, $ p_K^{(n+1)} \geq 0 $ and Theorem 8. is proved.

\subsection{Selection of the diffusion control parameters}
\label{subsec 7.3}

Based on the above constraints, we choose the diffusion control parameters as follows:
\begin{subequations} \label{definition_of_diffusion_control_parameters_for_variant_2}
    \begin{align}
        \left( \phi^{(u)}_e \right)_K &= \frac{1}{|(u_n)_K|} \left( \max \{ |u_n|_K, |u_n|_{K_e} \} + \frac{\sqrt{\rho_K} (u_n)_K + \sqrt{\rho_{K_e}} (u_n)_{K_e}}{\sqrt{\rho_K} + \sqrt{\rho_{K_e}}} \right), \label{2D_phi_1} \\
        \left( \phi^{(c)}_e \right)_K &= \exp \left( -\frac{(c_n)_K^2}{(u_n)_K^2} \right). \label{2D_phi_2}
    \end{align}
\end{subequations}
The details of the 2D numerical flux corresponding to the chosen $ \phi $'s are given in Appendix B. For the remaining choices of the diffusion control parameters in Equations \eqref{variant 1} and \eqref{variant 2}, the numerical results exhibit essentially the same behavior. Since no appreciable differences are observed in the resulting solutions, we restrict our presentation to the above choice.

\vspace{1\baselineskip}
\textbf{Theorem 9.} Let the states in the finite-volume cells adjacent to an interface satisfy the conditions for a stationary slip surface. Then the modified PVU Boltzmann scheme equipped with the diffusion control parameters \eqref{definition_of_diffusion_control_parameters_for_variant_2} admits this discontinuity as a discrete equilibrium. In particular,
\begin{equation} \label{2D discrete equilibrium}
    \bar{\boldsymbol{U}}_K^{(n + 1)} = \bar{\boldsymbol{U}}_K^{(n)}
\end{equation}
holds for all $ K $ and all $ n \geq 0 $.

\vspace{1\baselineskip}
Proof: We have, from \eqref{2D_finite_volume_update_formula},
\begin{equation}
    \bar{\boldsymbol{U}}_K^{n + 1} = \bar{\boldsymbol{U}}_K^{n} - \sigma \sum_e \hat{\textbf{f}}_e,
\end{equation}
where the interface numerical flux can be expressed as
\begin{equation} \label{e flux}
    \hat{\textbf{f}}_e = \dfrac{1}{2} \left( \boldsymbol{G}_{n, K} + \boldsymbol{G}_{n, K_e} \right) - \boldsymbol{d}_{e}^{(a)} - \boldsymbol{d}_{e}^{(p)}.
\end{equation}

From \eqref{2D_definition_advective_numerical_flux_PVU_plus}, we note that
\begin{equation}
    \boldsymbol{d}_{e}^{(a)} = \left( \phi_{e}^{(u)} \, | u_n | \, \boldsymbol{U} \right)_{K_e} - \left( \phi_{e}^{(u)} \, | u_n | \, \boldsymbol{U} \right)_K
\end{equation}

From \eqref{2D_phi_1}, we note that
\begin{equation}
    \lim_{u_n \rightarrow0} \phi^{(u)} \, | u_n | = 0.
\end{equation}
Thus,
\begin{equation}
    \boldsymbol{d}_{e}^{(a)} = 0.
\end{equation}
From \eqref{2D_phi_2}, we note that
\begin{equation}
    \lim_{u \rightarrow 0} \phi^{(c)} = 0,
\end{equation}
which implies that
\begin{equation}
    \boldsymbol{d}_{e}^{(p)} = 0.
\end{equation}

The averaged flux appearing in \eqref{e flux} reduce to
\begin{equation}
    \hat{\textbf{f}}_{e} = \begin{bmatrix}
        0 & p_0 & 0 & 0
    \end{bmatrix}^T.
\end{equation}
Therefore, the finite volume update reduces to \eqref{2D discrete equilibrium} and the modified PVU scheme preserves the stationary slip surface as a discrete equilibrium.


\section{Extension to second-order accuracy in two-dimensions}
\label{sec 8}

In this section, we use the MUSCL-like strategy to extend the modified PVU Boltzmann scheme based on hat distribution functions to second-order accuracy in space and time. In order to achieve second-order accuracy in space, a piecewise-linear approximation of the conserved variable vector is assumed in each cell, \textit{i.e.},
\begin{equation} \label{2D_2ndOrder_Reconstruction}
    \boldsymbol{U} (\boldsymbol{x}) = \boldsymbol{U}_K + \nabla \boldsymbol{U} \Big|_K \left( \boldsymbol{x} - \boldsymbol{x}_K \right),
\end{equation}
for all $ \boldsymbol{x} \in K $, $ \nabla ( \cdot ) $ denotes the gradient of the quantity in parenthesis and $ \boldsymbol{x}_K $ denote the coordinates of the cell centroid. In each cell, the gradient is evaluated using the approximate form of the exact relation,
\begin{equation}
    \int_{K} \, dA \, \nabla ( \cdot) = \oint_{\partial K} \, ds \, ( \cdot) \, \hat{\boldsymbol{n}}.
\end{equation}
Thus, the conserved variables are reconstructed at each cell interface $ e $ as
\begin{equation} \label{2D_2ndOrder_Reconstruction_interface}
\begin{aligned}
    \rho^L_{e} &= \rho_K + \nabla \rho \big|_{K} \bigl( \boldsymbol{x}_e - \boldsymbol{x}_K \bigr), \quad \quad \rho^R_{e} = \rho_{K_e} + \nabla \rho \big|_{K_e} \bigl( \boldsymbol{x}_e - \boldsymbol{x}_{K_e} \bigr), \\
    \boldsymbol{m}^L_{e} &= \boldsymbol{m}_K + \nabla \boldsymbol{m} \big|_{K} \bigl( \boldsymbol{x}_e - \boldsymbol{x}_K \bigr), \quad \boldsymbol{m}^R_{e} = \boldsymbol{m}_{K_e} + \nabla \boldsymbol{m} \big|_{K_e} \bigl( \boldsymbol{x}_e - \boldsymbol{x}_{K_e} \bigr), \\
    \mathcal{E}^L_{e} &= \mathcal{E}_{K} + \nabla \mathcal{E} \big|_{K} \bigl( \boldsymbol{x}_e - \boldsymbol{x}_K \bigr), \quad \quad \mathcal{E}^R_{e} = \mathcal{E}_{K_e} + \nabla \mathcal{E} \big|_{K_e} \bigl( \boldsymbol{x}_e - \boldsymbol{x}_{K_e} \bigr),
\end{aligned}
\end{equation}
where $ \boldsymbol{m}, \mathcal{E} $ denote the momentum vector and the total energy, respectively, while $ ( \cdot )^L_{e}, \, ( \cdot)^R_{e} $ denote the conserved variables reconstructed on the left and right of the interface, respectively. A semi-discrete finite volume discretization of the two-dimensional collisionless Boltzmann-like equations reads

\begin{subequations} \label{eq:2D_semi_discrete_update}
    \begin{align}
        \int_{A} \, \frac{\partial f}{\partial t} \, dA &= -\sum_{e} \left( \boldsymbol{u}_e \cdot \hat{\boldsymbol{n}}_e  + \boldsymbol{c}_e \cdot \hat{\boldsymbol{n}}_e \right) \, f_e \, m(e), \\
        \int_{A} \, \frac{\partial g}{\partial t} \, dA &= -\sum_{e} \left( \boldsymbol{u}_e \cdot \hat{\boldsymbol{n}}_e  + \boldsymbol{c}_e \cdot \hat{\boldsymbol{n}}_e \right) \, g_e \, m(e).
    \end{align}
\end{subequations}
The hat function is constructed at each interface using the reconstructed conserved variables, which read
\begin{subequations}
    \begin{align}
        f^L_{e} &= \bar{f}^{eq} (\rho^L_{e}, \, \boldsymbol{u}^L_{e}, \, p^L_{e}), \quad g^L_{e} = \bar{g}^{eq} (\rho^L_{e}, \, \boldsymbol{u}^L_{e}, \, p^L_{e}), \\
        f^R_{e} &= \bar{f}^{eq} (\rho^R_{e}, \, \boldsymbol{u}^R_{e}, \, p^R_{e}), \quad g^R_{e} = \bar{g}^{eq} (\rho^R_{e}, \, \boldsymbol{u}^R_{e}, \, p^R_{e}).
    \end{align}
\end{subequations}
The kinetic numerical flux normal to each interface is defined based on the reconstructed hat function, \textit{i.e.},
\begin{multline} \label{numerical_flux_1}
    (\boldsymbol{u}_e \cdot \hat{\boldsymbol{n}}_e) \, f_e = \left( \dfrac{(\boldsymbol{u}_e^{L} \cdot \hat{\boldsymbol{n}}_e) + |\boldsymbol{u}_e^{L} \cdot \hat{\boldsymbol{n}}_e| \, \phi_e^{(u)}}{2} \right) f^L_{e} \\
    +  \left( \dfrac{(\boldsymbol{u}_e^{R} \cdot \hat{\boldsymbol{n}}_e) - |\boldsymbol{u}_e^{R} \cdot \hat{\boldsymbol{n}}_e| \, \phi_e^{(u)}}{2} \right) f^R_{e},
\end{multline}
\begin{multline} \label{numerical_flux_2}
    (\boldsymbol{u}_e \cdot \hat{\boldsymbol{n}}_e) \, g_e = \left( \dfrac{(\boldsymbol{u}_e^{L} \cdot \hat{\boldsymbol{n}}_e) + |\boldsymbol{u}_e^{L} \cdot \hat{\boldsymbol{n}}_e| \, \phi_e^{(u)}}{2} \right) g^L_{e} \\
    +  \left( \dfrac{(\boldsymbol{u}_e^{R} \cdot \hat{\boldsymbol{n}}_e) - |\boldsymbol{u}_e^{R} \cdot \hat{\boldsymbol{n}}_e| \, \phi_e^{(u)}}{2} \right) g^R_{e},
\end{multline}
\begin{multline} \label{numerical_flux_3}
     (\boldsymbol{c}_e \cdot \hat{\boldsymbol{n}}_e) \, f_e = \left( \dfrac{(\boldsymbol{c}_e^{L} \cdot \hat{\boldsymbol{n}}_e) + |\boldsymbol{c}_e^{L} \cdot \hat{\boldsymbol{n}}_e| \, \phi_e^{(c)}}{2} \right) f^L_{e} \\
    +  \left( \dfrac{(\boldsymbol{c}_e^{R} \cdot \hat{\boldsymbol{n}}_e) - |\boldsymbol{c}_e^{R} \cdot \hat{\boldsymbol{n}}_e| \, \phi_e^{(c)}}{2} \right) f^R_{e},
\end{multline}
\begin{multline} \label{numerical_flux_4}
     (\boldsymbol{c}_e \cdot \hat{\boldsymbol{n}}_e) \, g_e = \left( \dfrac{(\boldsymbol{c}_e^{L} \cdot \hat{\boldsymbol{n}}_e) + |\boldsymbol{c}_e^{L} \cdot \hat{\boldsymbol{n}}_e| \, \phi_e^{(c)}}{2} \right) g^L_{e} \\
    +  \left( \dfrac{(\boldsymbol{c}_e^{R} \cdot \hat{\boldsymbol{n}}_e) - |\boldsymbol{c}_e^{R} \cdot \hat{\boldsymbol{n}}_e| \, \phi_e^{(c)}}{2} \right) g^R_{e}.
\end{multline}
Linear reconstruction of the conserved variable vector within a finite volume cell may produce unphysical values, such as negative density and pressure, at cell interfaces. These unphysical states, in turn, violate the positivity of the equilibrium distribution function, which depends on density $ \rho $, velocity $ \boldsymbol{u} $ and pressure $ p $. To address this issue and ensure the positivity of density and pressure at interfaces, we adopt a two-stage limiting strategy: the Venkatakrishnan limiter is applied first to reduce spurious oscillations, followed by the Zhang-Shu limiter, which enforces positivity preservation. We refer the reader to \cite{Venkatakrishnan} and \cite{Zhang_Shu} for the implementation details of the respective limiters on rectangular grids. With second-order accuracy in the spatial dimension, we choose the three-stage Strong Stability Preserving Runge-Kutta (SSPRK-3) method for explicit time discretization of the temporal term in Equation (\ref{eq:2D_semi_discrete_update}) which reads
\begin{subequations} \label{SSPRK-3_update_formula}
    \begin{equation} \label{Euler_forward_in_time}
        h^{(1)-}_K = \bar{h}^{eq}_K - \frac{\Delta t}{m(K)} \, R_K (\bar{h}^{eq})
    \end{equation}
    \begin{equation}
        h^{(2)-}_K = \frac{3}{4} \bar{h}^{eq}_K + \frac{1}{4} h^{(1)}_K - \frac{1}{4} \frac{\Delta t}{m(K)} \, R_K (h^{(1)})
    \end{equation}
    \begin{equation}
        h^{(n+1)-}_K = \frac{1}{3} \bar{h}^{eq}_K + \frac{2}{3} h^{(2)}_K -\frac{2}{3} \frac{\Delta t}{m(K)} \, R_K (h^{(2)}),
    \end{equation}
\end{subequations}
for $ h = f $ or $ g $, and $ R $ denotes the summation of the numerical fluxes at all the edges $ e $ of the cell $ K $, termed residual. Using Equation (\ref{SSPRK-3_update_formula}) in (\ref{2D_definition_of_conserved_variable}) leads to the macroscopic update of the conserved variable vector, \textit{i.e.},
\begin{subequations}
    \begin{equation}
        \bar{\boldsymbol{U}}^{(1)}_K = \bar{\boldsymbol{U}}^{n}_K - \frac{\Delta t}{m(K)} \boldsymbol{R}_K (\bar{\boldsymbol{U}}^n)
    \end{equation}
    \begin{equation}
        \bar{\boldsymbol{U}}^{(2)}_K = \frac{3}{4} \bar{\boldsymbol{U}}^n_K + \frac{1}{4} \bar{\boldsymbol{U}}^{(1)}_K - \frac{1}{4} \frac{\Delta t}{m(K)} \boldsymbol{R}_K (\bar{\boldsymbol{U}}^{(1)})
    \end{equation}
    \begin{equation}
        \bar{\boldsymbol{U}}^{(n+1)}_K = \frac{1}{3} \bar{\boldsymbol{U}}^n_K + \frac{2}{3} \bar{\boldsymbol{U}}^{(2)}_K - \frac{2}{3} \frac{\Delta t}{m(K)} \boldsymbol{R}_K (\bar{\boldsymbol{U}}^{(2)}).
    \end{equation}
\end{subequations}

\vskip 1em
\noindent
Note: It is to be noted that in each cell, the collision step instantaneously re-establishes the equilibrium without changing the mass, momentum, and energy densities after each Runge-Kutta stage for the streaming step.

\subsection{Positivity analysis of the second-order scheme}
\label{subsec 8.1}
\noindent
\textbf{Theorem 10.} Assume that $ \Bar{f}^{eq}_K = 0 $ for $ |(c_n)_K| \geq |(c_n)_K|_{max} $. Then, under the CFL-like condition,
\begin{equation} \label{2D_CFL_like_condition_second_order}
    \frac{\Delta t}{m(K)} \max_{e} \{ | \boldsymbol{u}_e \cdot \hat{\boldsymbol{n}}_e | + | \boldsymbol{c}_e \cdot \hat{\boldsymbol{n}}_e |_{max} \} \sum_{e} m(e) \leq \dfrac{1}{4},
\end{equation}
the Euler forward-in-time, second-order-in-space modified PVU Boltzmann scheme (\ref{Euler_forward_in_time}) preserves positivity, \textit{i.e.}, 
\begin{equation}
    \rho^{(n + 1)} > 0, \; p^{(n + 1)} > 0 \quad \text{if} \quad \rho^{(n)} > 0, \; p^{(n)} > 0. 
\end{equation}

\noindent
Proof: We have, from (\ref{Euler_forward_in_time}), (\ref{numerical_flux_1}), and (\ref{numerical_flux_3}),
\begin{multline} \label{second_order_update_formula}
    f_K^{(n + 1)-} = \Bar{f}_K^{eq} - \frac{\Delta t}{2 m(K)} \sum_{e} \left( \boldsymbol{u}_e^L \cdot \hat{\boldsymbol{n}}_e + | \boldsymbol{u}_e^L \cdot \hat{\boldsymbol{n}}_e | \, \phi^{(u)}_e \right) \, f_e^L \, m(e) \\
    - \frac{\Delta t}{2 m(K)} \sum_{e} \left( \boldsymbol{c}_e^L \cdot \hat{\boldsymbol{n}}_e + | \boldsymbol{c}_e^L \cdot \hat{\boldsymbol{n}}_e | \, \phi^{(c)}_e \right) f_e^L \, m(e) - \frac{\Delta t}{2 m(K)} \sum_{e} h_{e}^- m(e),
\end{multline}
where
\begin{equation}
    h_{e}^- = \left( \boldsymbol{u}_e^R \cdot \hat{\boldsymbol{n}}_e - | \boldsymbol{u}_e^R \cdot \hat{\boldsymbol{n}}_e | \, \phi^{(u)}_e \right) f^{R}_{e} + \left( \boldsymbol{c}_{e}^R \cdot \hat{\boldsymbol{n}}_e - | \boldsymbol{c}_e^R \cdot \hat{\boldsymbol{n}}_e | \, \phi^{(c)}_e \right) f^{R}_{e}.
\end{equation}
We then note from Proposition 2. and (\ref{positivity_of_u_tilda_plus}), (\ref{positivity_of_u_tilda_minus}) that
\begin{equation}
    -\left\langle \frac{\Delta t}{2 m(K)} \sum_{e} h_{e}^- m(e) \right\rangle \geq 0.
\end{equation}
We note, at least for rectangular grids, that
\begin{equation} \label{cell_average_from_reconstructions}
    \left \langle \bar{f}^{eq}_K \right \rangle = \dfrac{1}{4} \left \langle \sum_e \, f_e^{L}  \right \rangle.
\end{equation}
Using (\ref{cell_average_from_reconstructions}) in (\ref{second_order_update_formula}), we get
\begin{multline} \label{in_between_equation_177}
    f_K^{(n + 1)-} = \dfrac{1}{4} \sum_e f_e^L - \frac{\Delta t}{2 m(K)} \sum_{e} \left( \boldsymbol{u}_e^L \cdot \hat{\boldsymbol{n}}_e + | \boldsymbol{u}_e^L \cdot \hat{\boldsymbol{n}}_e | \, \phi^{(u)}_e \right) \, f_e^L \, m(e) \\
    - \frac{\Delta t}{2 m(K)} \sum_{e} \left( \boldsymbol{c}_e^L \cdot \hat{\boldsymbol{n}}_e + | \boldsymbol{c}_e^L \cdot \hat{\boldsymbol{n}}_e | \, \phi^{(c)}_e \right) f_e^L \, m(e) - \frac{\Delta t}{2 m(K)} \sum_{e} h_{e}^- m(e).
\end{multline}
From (\ref{constraint_on_phi2}), we have
\begin{equation} \label{restriction_on_peculiar_velocity}
    \dfrac{1}{2} \left( \boldsymbol{c}_e^L \cdot \hat{\boldsymbol{n}}_e + | \boldsymbol{c}_e^L \cdot \hat{\boldsymbol{n}}_e | \, \phi_e^{(c)} \right) \leq | \boldsymbol{c}_e^L \cdot \hat{\boldsymbol{n}}_e | \leq \max_e \left\{ | \boldsymbol{c}_e \cdot \hat{\boldsymbol{n}}_e |_{max} \right\}.
\end{equation}
For a specific choice of $ \phi^{(1)} $ given in (\ref{2D_phi_1}), we note
\begin{equation} \label{restriction_on_bulk_velocity}
    \dfrac{1}{2} \left( \boldsymbol{u}_e^L \cdot \hat{\boldsymbol{n}}_e + | \boldsymbol{u}_e^L \cdot \hat{\boldsymbol{n}}_e | \, \phi_e^{(u)} \right) \leq \max_e \left\{ | \boldsymbol{u}_e \cdot \hat{\boldsymbol{n}}_e | \right\}.
\end{equation}
Multiplying equations (\ref{restriction_on_peculiar_velocity}) and (\ref{restriction_on_bulk_velocity}) by $ m(e) $ and $ f_e^L $, and summing over all edges $ e $, we obtain:
\begin{multline} \label{in_between_equation_180}
    \dfrac{1}{2} \sum_e \, \left( \left( \boldsymbol{u}_e^L \cdot \hat{\boldsymbol{n}}_e + | \boldsymbol{u}_e^L \cdot \hat{\boldsymbol{n}}_e | \, \phi_e^{(u)} \right) + \left( \boldsymbol{c}_e^L \cdot \hat{\boldsymbol{n}}_e + | \boldsymbol{c}_e^L \cdot \hat{\boldsymbol{n}}_e | \, \phi_e^{(c)} \right) \right) f_e^L \, m(e) \\
    \leq \left( \max_e \left\{ | \boldsymbol{c}_e \cdot \hat{\boldsymbol{n}}_e |_{max} \right\} + \max_e \left\{ | \boldsymbol{u}_e \cdot \hat{\boldsymbol{n}}_e | \right\} \right) \sum_e \, f_e^L \, m(e).
\end{multline}
Using Equation (\ref{in_between_equation_180}) in (\ref{in_between_equation_177}), we note that the following condition must be satisfied to ensure the positivity of the first term:
\begin{equation}
    \dfrac{\Delta t}{m(K)} \max_e \left\{ | \boldsymbol{u}_e \cdot \hat{\boldsymbol{n}}_e | + | \boldsymbol{c}_e \cdot \hat{\boldsymbol{n}}_e |_{max} \right\} \, m(e) \leq \dfrac{1}{4} \quad \forall e.
\end{equation}
Thus, by a CFL-like condition given in Equation (\ref{2D_CFL_like_condition_second_order}), we have $ f_K^{(n + 1)-} \geq 0 $ and hence
\begin{subequations}
    \begin{align}
        \rho^{(n + 1)}_K &= \int_{\boldsymbol{c}} \, D \boldsymbol{c} \, f_K^{(n + 1)-} \geq 0, \\
        p_K^{(n + 1)} &= \int_{\boldsymbol{c}} \, D \boldsymbol{c} \, c_K^2 \, f_K^{(n + 1)-} \geq 0.
    \end{align}
\end{subequations}


\section{Results}
\label{sec 9}

\subsection{Experimental Order of Convergence}
\label{subsec 9.1}

The experimental order of convergence (EOC) provides a practical means of verifying the convergence order of a scheme. For this purpose a smooth test is considered on the interval $ x \in [0, 1] $, where density varies sinusoidally while velocity and pressure remain constant. The analytical solution of this problem reads
\begin{equation} \label{exact_solution_of_density_for_EOC_test}
    \rho(x,t) = \rho_0 (x - u t) = 1.0 + 0.2 \sin{\{ 2 \pi (x - 0.1 t) \}},
\end{equation}
\begin{equation} \label{exact_solution_of_velocity_pressure_for_EOC_test}
    u(x, t) = 0.1, \, p(x, t) = 0.5.
\end{equation}
The numerical solution is allowed to evolve until t = 0.5 units for the grid of sizes $ N_x $ = 10, 20, 40, 80, 160, 320, 640, 1280, with periodic boundary conditions at the left and right ends of the one-dimensional domain. The $ L_1 $ and $ L_2 $ norms of the error are evaluated as
\begin{equation}
    \bigl|\bigl| \varepsilon \bigl|\bigl|_{L_1} = {\dfrac{1}{N_x}} \sum^{N}_{i=1} \bigl| \rho_i - \rho^{exact}_i \bigl|,
\end{equation}
\begin{equation}
    \bigl|\bigl| \varepsilon\bigl |\bigl|_{L_2} = \sqrt{{\dfrac{1}{N_x}} \sum^{N}_{i=1} \bigl( \rho_i - \rho^{exact}_i \bigl)^2},
\end{equation}
where $ \rho_i $ and $ \rho^{exact}_i $ are the numerical and exact densities at the cell center $ i $. Pressure and velocity remain constant, with no contribution to the error.
The experimental order of convergence is then evaluated as
\begin{equation}
    EOC = log_2 \Biggl( \dfrac{\bigl|\bigl| {\varepsilon}_{{N_x}/2} \bigl|\bigl|}{\bigl|\bigl| {\varepsilon}_{N_x} \bigl|\bigl|} \Biggl).
\end{equation}
The error norms $ L_1 $ and $ L_2 $ and the EOC evaluated for the respective grid sizes are listed in Table \ref{tab:1D_L1_L2_variant2}.

\begin{table}[h]
    \centering
    \setlength{\tabcolsep}{10pt}
    \renewcommand{\arraystretch}{1.25}
    \begin{tabular}{ l l c c c c }
    \hline
        $ N_x $ & $ dx $ & $ L_1 $ error & EOC & $ L_2 $ error & EOC \\
    \hline
        10 & 0.1  & 0.0015 & & 0.0018 & \\
        20 & 0.05 & 3.46E-04 & 2.11 & 3.87E-04 & 2.22 \\
        40 & 0.025 & 8.34E-05 & 2.05 & 9.27E-05 & 2.06 \\
        80 & 0.0125 & 2.06E-05 & 2.01 & 2.29E-05 & 2.02 \\
        160 & 0.00625 & 5.14E-06 & 2.00 & 5.72E-06 & 2.00 \\
        320 & 0.003125 & 1.29E-06 & 2.00 & 1.43E-06 & 2.00 \\
        640 & 0.0015625 & 3.21E-07 & 2.00 & 3.57E-07 & 2.00 \\
        1280 & 0.00078125 & 8.03E-08 & 2.00 & 8.92E-08 & 2.00 \\
    \hline
    \end{tabular}
    \caption{EOC study for the second-order modified PVU Boltzmann scheme; CFL = $ 0.2 $}
    \label{tab:1D_L1_L2_variant2}
\end{table}

\subsection{1D results and discussions}
\label{subsec 9.2}

In this subsection, we evaluate performance of the modified PVU Boltzmann scheme on several one-dimensional Riemann problems for which analytical solutions are available. Each test consists of constant left and right states ($ \textbf{U}_L, \textbf{U}_R $) separated at the point $ x_0 $. Computations are performed on the domain $ x \in [0, 1] $ discretized with $ 200 $ cells. CFL numbers of 0.4 and 0.2 are used to obtain first-order and second-order results, respectively. The initial conditions for all numerical tests are cataloged in the Table \ref{tab:1D_Riemann_problems}.

\renewcommand{\arraystretch}{1.5}
\begin{table}[H]
    \centering
    \begin{tabular}{|c|c|c|c|c|c|c|c|c|}
        \hline
        Test & $ \rho_L $ & $ u_L $ & $ p_L $ & $ \rho_R $ & $ u_R $ & $ p_R $ & runtime & $ x_0 $ \\
        \hline
        1 & 1.4 & 0.0 & 1.0 & 1.0 & 0.0 & 1.0 & 2.0 & 0.5 \\
        \hline
        2 & 1.4 & 0.1 & 1.0 & 1.0 & 0.1 & 1.0 & 1.0 & 0.5 \\
        \hline
        3 & 3.86 & -0.81 & 10.33 & 1.0 & -3.44 & 1.0 & 4.0 & 0.1 \\
        \hline
        4 & 1.0 & 0.75 & 1.0 & 0.125 & 0.0 & 0.1 & 0.2 & 0.3 \\
        \hline
        5 & 1.0 & -2.0 & 0.4 & 1.0 & 2.0 & 0.4 & 0.15 & 0.5 \\
        \hline
        6 & 5.999 & 19.5975 & 460.894 & 5.992 & -6.1963 & 46.095 & 0.035 & 0.4 \\
        \hline
    \end{tabular}
    \caption{1D Riemann problems}
    \label{tab:1D_Riemann_problems}
\end{table}

\vspace{-0.5 cm}
\begin{figure}[H]
    \centering
    \includegraphics[width=1\linewidth]{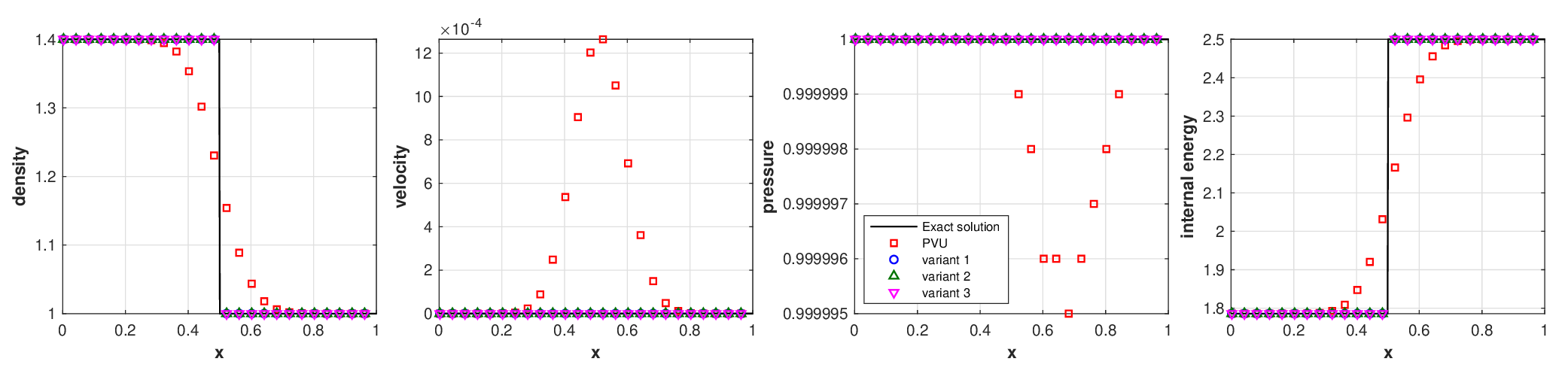}
    \vspace{-0.8 cm}
    \caption{Stationary contact discontinuity}
    \label{1D_figure1}
\end{figure}

Numerical test 1 evaluates the performance of a scheme on an isolated, grid-aligned stationary contact discontinuity. As shown in Figure \ref{1D_figure1}, all the proposed variants of the modified PVU Boltzmann scheme preserve the stationary contact discontinuity as a discrete equilibrium, thereby eliminating spurious pressure and velocity oscillations observed with the PVU scheme.

\begin{figure}[H]
    \centering
    \includegraphics[width=1\linewidth]{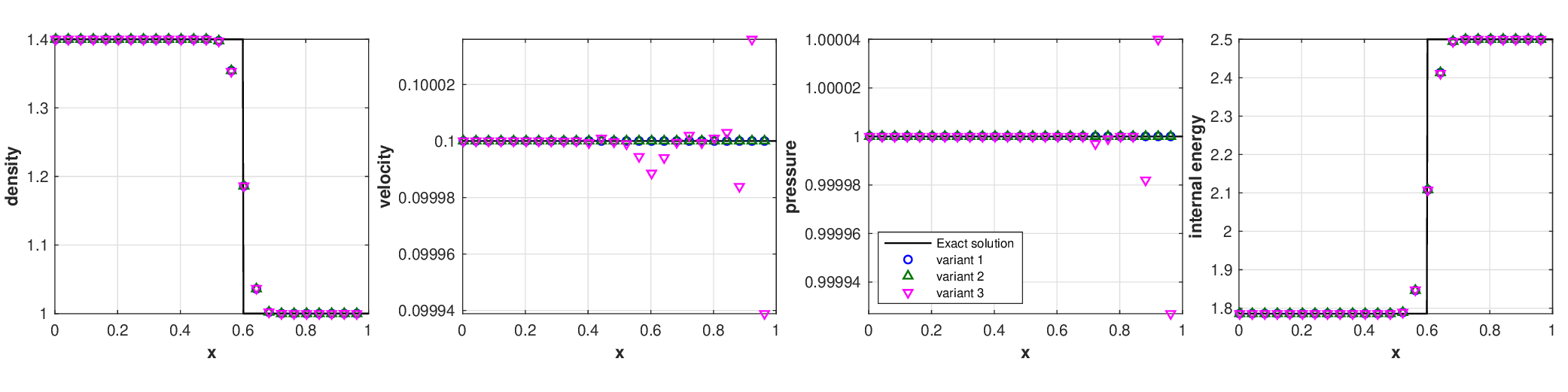}
    \vspace{-0.8 cm}
    \caption{Slowly moving contact discontinuity}
    \label{1D_figure2}
\end{figure}

Numerical test 2 assesses the performance of a scheme across an isolated, slowly moving contact discontinuity. The solution consists of a right-moving contact discontinuity propagating at a constant speed. As shown in Figure \ref{1D_figure2}, the results corresponding to the choice of $ \phi^{(c)} $ in \eqref{variant 3} exhibit velocity and pressure oscillations of $\mathcal{O}\left(10^{-5}\right)$. In contrast, the results obtained using the alternative choices of $ \phi^{(c)} $ in \eqref{variant 1} and \eqref{variant 2} remain free of such oscillations. This test therefore represents the only case among those considered in which the choice of the diffusion control parameter has a discernible effect on the numerical solution. Consequently, all subsequent one-dimensional results are reported for variant 1 of the modified PVU Boltzmann scheme.

\begin{figure}[H]
    \centering
    \includegraphics[width=1\linewidth]{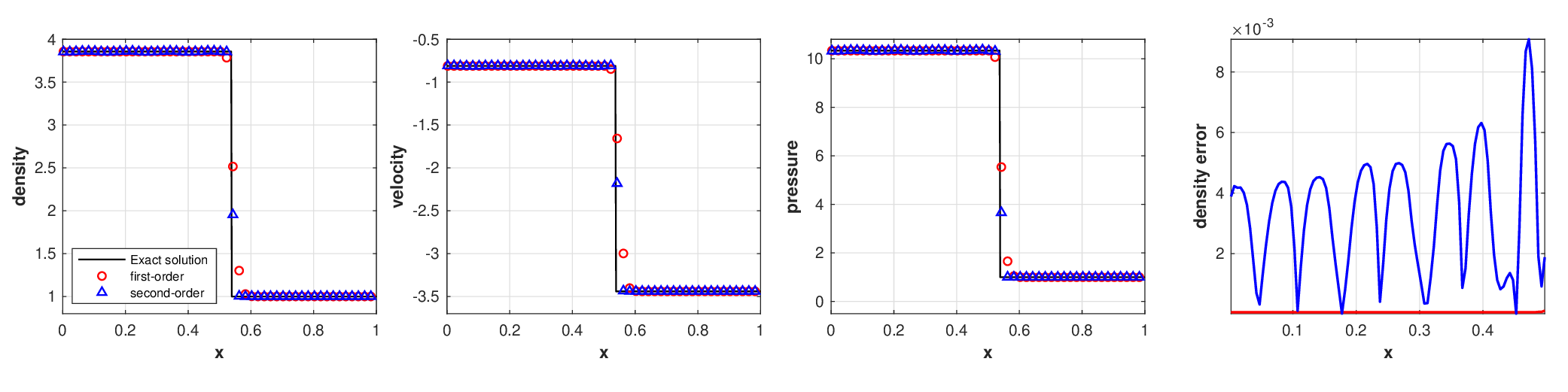}
    \vspace{-0.8 cm}
    \caption{Slowly moving shock}
    \label{1D_figure3}
\end{figure}

The solution of Test 3 consists of a slowly moving shock propagating to the right. This test is designed to assess the scheme's ability to avoid low-frequency post-shock oscillations. Figure \ref{1D_figure3} compares the first- and second-order results obtained with the modified PVU Boltzmann scheme. The first-order results do not exhibit any post-shock oscillations, demonstrating the scheme's ability to capture the shock without introducing spurious fluctuations. In contrast, the second-order results exhibit subtle post-shock oscillations of the order $ \mathcal{O} (10^{-3}) $.

\begin{figure}[H]
    \centering
    \includegraphics[width=1\linewidth]{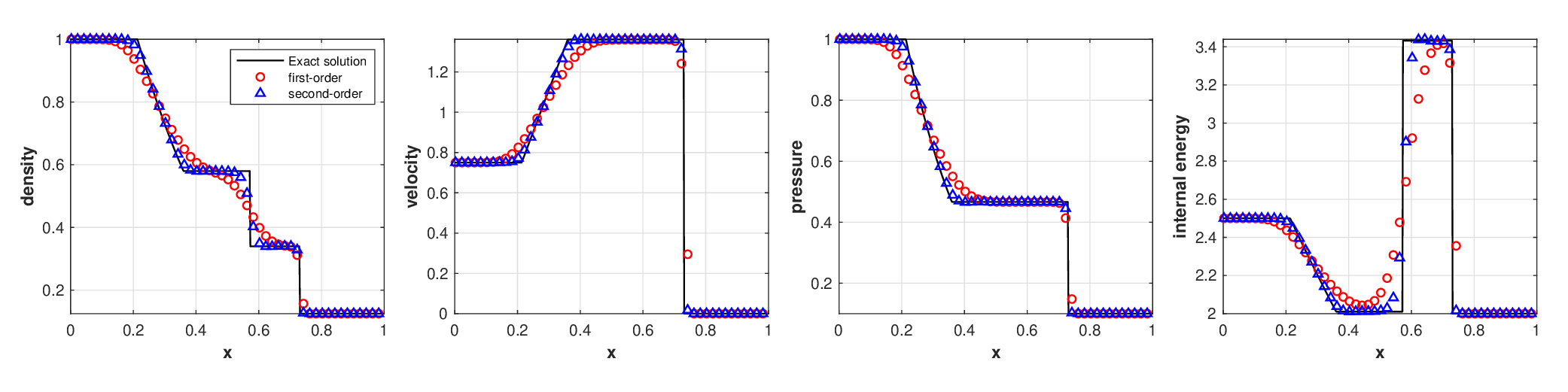}
    \vspace{-1 cm}
    \caption{Sod's shock tube problem}
    \label{1D_figure4}
\end{figure}

Numerical test 4 is the modified version of the Sod's shock tube problem; its analytical solution consists of a left-sonic expansion, a right-moving contact discontinuity, and a right-moving shock. This test assesses the entropy satisfaction property of a numerical method. As can be seen in Figure \ref{1D_figure4}, neither the first- nor the second-order results of the modified PVU Boltzmann scheme exhibit entropy violating solutions.

\begin{figure}[H]
    \centering
    \includegraphics[width=1\linewidth]{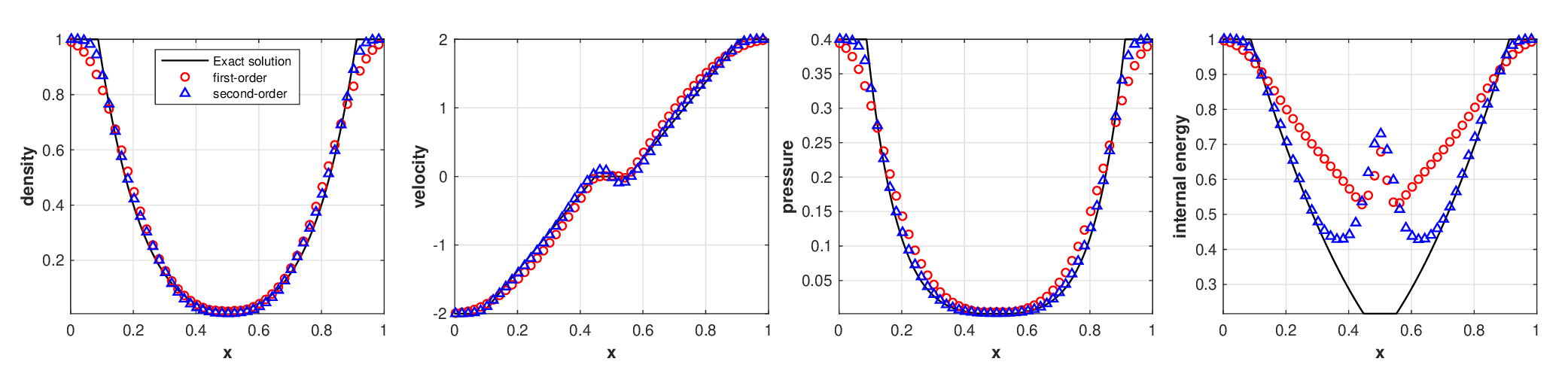}
    \vspace{-0.8 cm}
    \caption{Overheating problem}
    \label{1D_figure5}
\end{figure}

The exact solution of test 5 consists of left and right expansion fans connected by a stationary contact discontinuity. This test helps assess the performance of a numerical scheme near low-density regions. As seen in Figure \ref{1D_figure5}, the modified PVU Boltzmann scheme remains robust in the presence of low-density regions and does not exhibit numerical failure. However, an overshoot in the internal energy is observed near the center, a feature commonly encountered in upwind schemes.

\begin{figure}[H]
    \centering
    \includegraphics[width=1\linewidth]{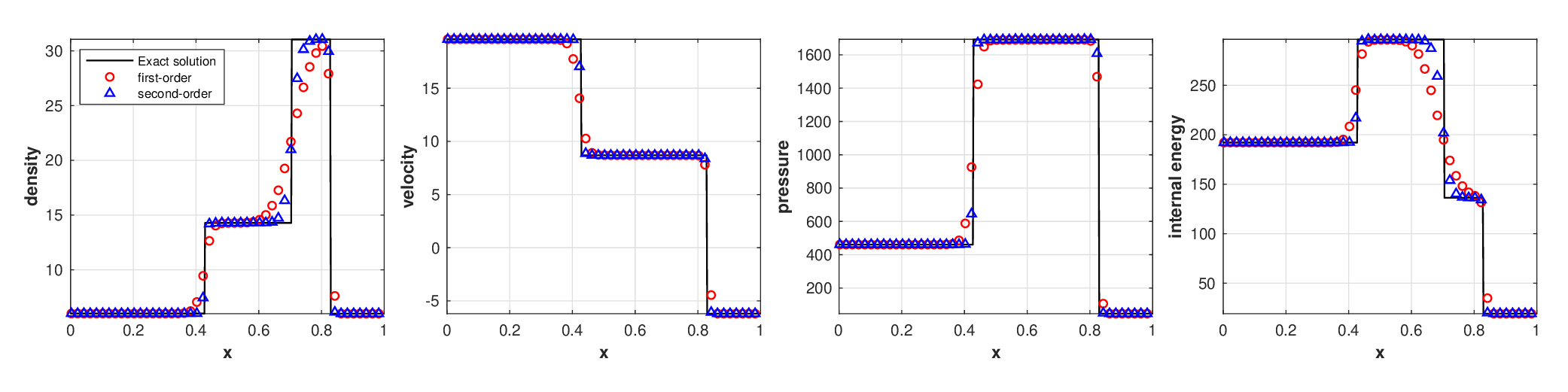}
    \vspace{-1 cm}
    \caption{Colliding strong shocks}
    \label{1D_figure6}
\end{figure}

Numerical test 6 is also designed to test the robustness of the numerical scheme. The solution consists of three strong discontinuities moving to the right. The modified PVU scheme captures the key flow features reasonably well at first- and second-order accuracy.

\subsection{2D results and discussions}
\label{subsec 9.3}
In this section, we test the (variant 3 of the) modified PVU scheme on standard 2D inviscid gas dynamics problems to assess its robustness, accuracy, and positivity-preserving properties. We use the kinetic characteristic boundary condition (KCBC) at the walls to ensure flow tangency. Unless otherwise mentioned, the CFL numbers of 0.4 and 0.2 are used to generate first-order and second-order results, respectively. For steady-state convergence, the first-order solution is advanced until the density residual reaches at least $ 10^{-8} $ is achieved, while the second-order solution is run until the residual reaches at least $ 10^{-6} $ is achieved.  We set K (a free parameter in the Venkatakrishnan limiter) to $ 0.3 $ for unsteady problems and $ 5 $ for steady problems, unless otherwise mentioned.

\subsubsection{Oblique shock reflection}
\label{subsubsec 9.3.1}

\begin{figure}[ht]
   \centering
   \includegraphics[scale=0.3]{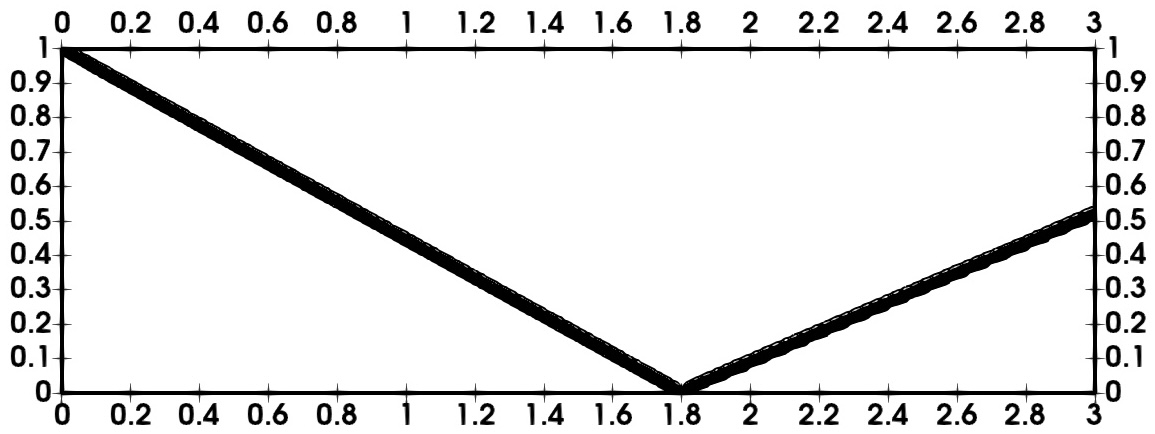}
   \caption{Oblique shock reflection; density contours (1:0.053:2.6)}
   \label{2D_figure1}
\end{figure}

This test considers an incoming oblique shock that strikes and reflects from a solid wall. The free-stream Mach number is 2.9. The computational domain is a rectangular region defined by $ [3, 0] \times [0, 1] $, discretized using $ 240 \times 80 $ Cartesian cells. A supersonic inlet boundary condition is applied at the left boundary, while a kinetic characteristic boundary condition (KCBC) is imposed along the bottom wall. Post-shock conditions, obtained from inviscid gas dynamics relations, are prescribed at the top boundary, and a zero gradient boundary condition is used at the right boundary. The domain is initialized with free-stream conditions. The density contours obtained using the second-order accurate modified PVU Boltzmann scheme are shown in Figure \ref{2D_figure1}. The oblique shocks are sharply captured.

\begin{figure}[ht]
   \centering
   \includegraphics[width=0.7\textwidth]{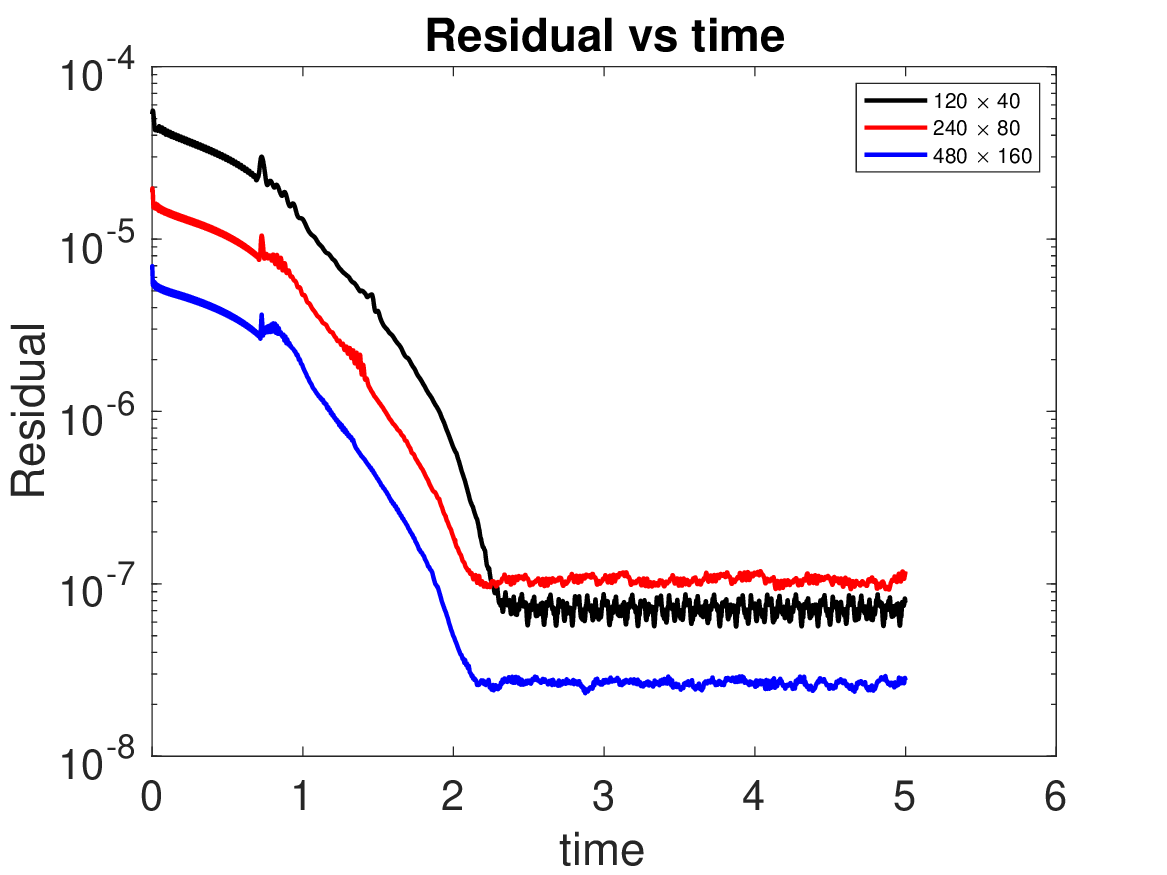}
   \caption{Oblique shock reflection test case: RMS error in density vs time}
    \label{2D_figure2}
\end{figure}
The convergence behavior of the proposed second-order scheme is assessed for the oblique shock reflection test case on three successively refined meshes: $ 120 \times 40 $, $ 240 \times 80 $ and $ 480 \times 160 $. Figure \ref{2D_figure2} shows the residuals as functions of time, representing the time-marching iterations for second-order simulations. With mesh refinement, the residuals exhibit an overall reduction, although not strictly monotonic. In particular, a characteristic non-monotonic bump appears around time $ 0.8 $, which is typical for many schemes. The second-order solutions demonstrate rapid residual decay and reach comparably low residual levels across all grids, at least $ \mathcal{O}(10^{-7}) $, indicating good convergence properties and confirming the expected accuracy and robustness of the scheme.

\subsubsection{Horizontal slip flow}
\label{subsubsec 9.3.2}

\begin{figure}[ht]
   \centering
   \includegraphics[scale=0.2]{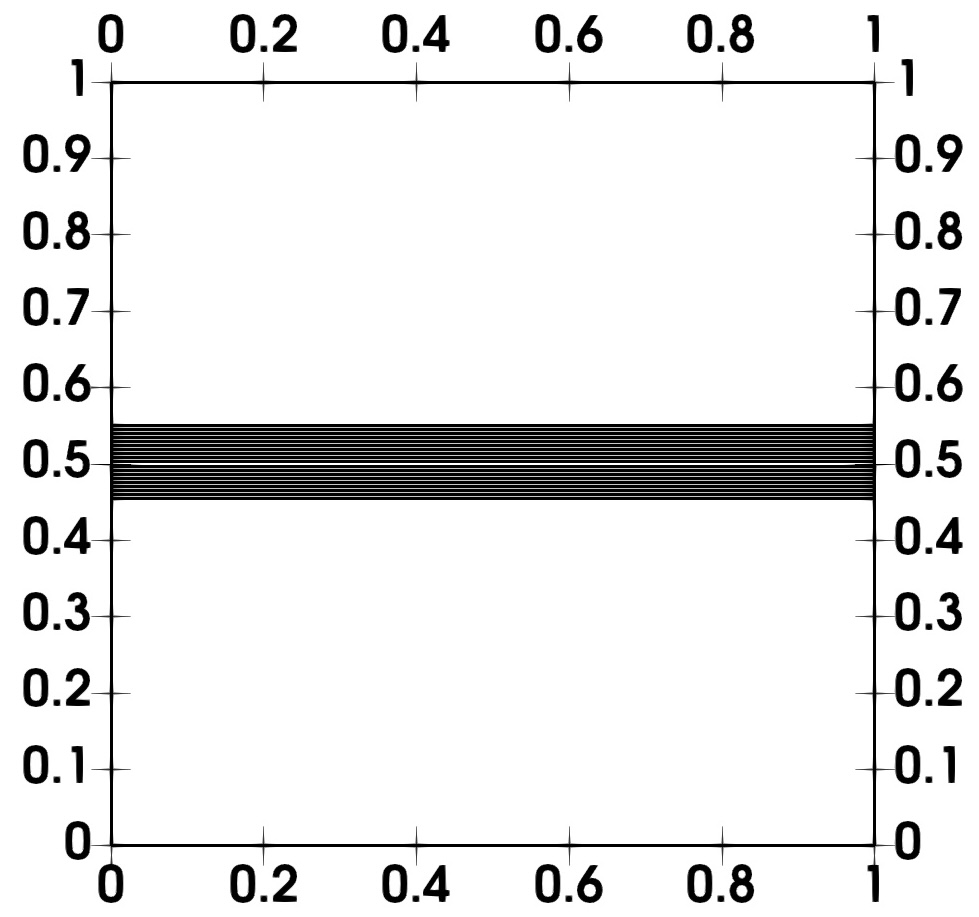}
   \caption{Horizontal slip flow; Horizontal velocity contours (2.0:0.05:3.0)}
   \label{2D_figure3}
\end{figure}
This test case involves a Mach $ 3.0 $ stream sliding past a Mach $ 2.0 $ stream. The computational region consists of a unit square whose extents are given by $ [1, 0] \times [0, 1] $. The computational domain is divided into $ 20 \times 20 $ finite-volume cells. At the initial instant the upper half of the domain $ y > 0.5 $ is filled with a horizontal Mach $ 3.0 $ flow while the lower half $ y < 0.5 $ contains a horizontal Mach $ 2.0 $ flow; density and pressure are held uniform throughout. The supersonic inlet boundary condition is used at the left and right boundaries. The zero gradient boundary condition is employed at the top and bottom boundaries. Figure \ref{2D_figure3} shows the contours of the horizontal velocity computed using the first-order modified PVU Boltzmann scheme. The result demonstrates that the proposed kinetic scheme accurately captures the grid-aligned, isolated, and steady slip surface without introducing numerical diffusion.

\subsubsection{Positivity test case}
\label{subsubsec 9.3.3}

\begin{figure}[H]
   \centering
   \includegraphics[scale=0.22]{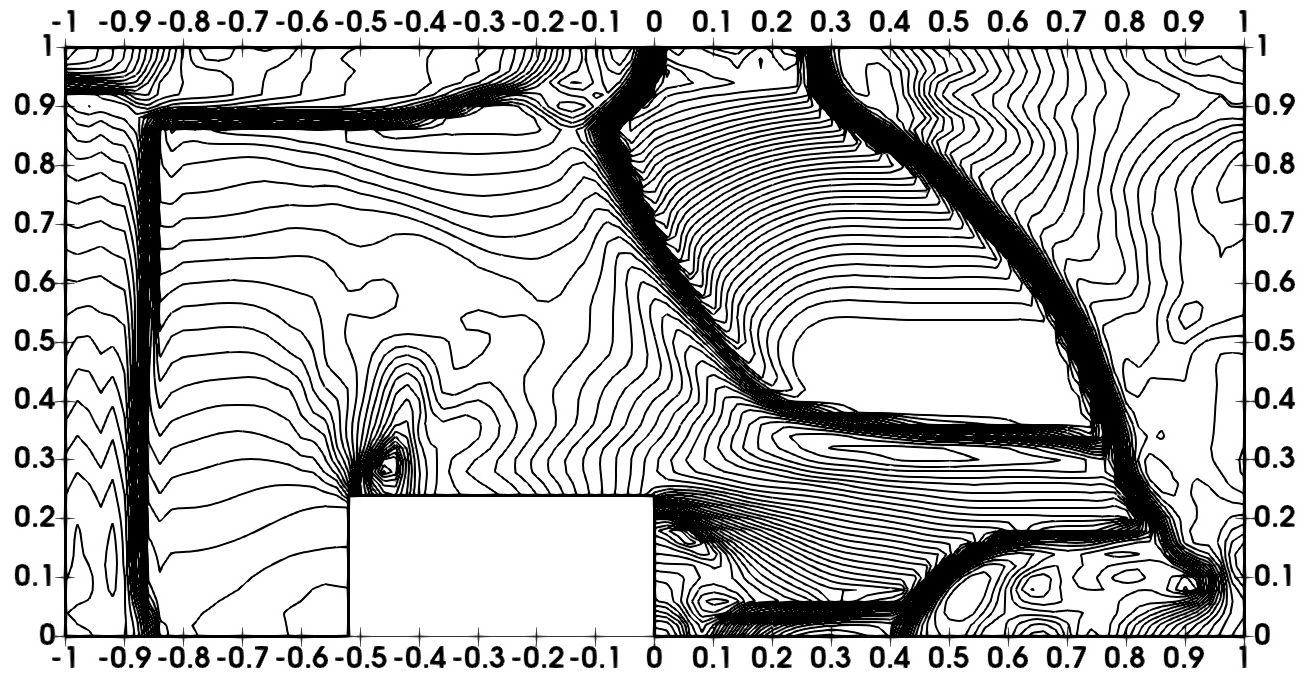}
   \caption{2D Enclosure test; 75 logarithmic pressure contours ($ 23.9079 $ - $ 3.25518 \times 10^6 $)}
   \label{2D_figure5}
\end{figure}
This test case is taken from Parent \cite{Parent} and is designed to assess the positivity-preserving property of a numerical scheme. This test consists of large initial Mach numbers, which leads to strong expansions in the first few iterations. It is an excellent test to empirically determine the positivity-preserving property of both variants of the modified PVU scheme. The flow is enclosed within $ x \in [-1, 1] $ and $ y \in [0, 1] $, with a cutout at $ x = -0.52 $ and a height of $ 0.24 $ units. The computational domain is divided into Cartesian cells, each with a size of $ 0.02 $ units along both dimensions. The kinetic characteristic boundary condition is employed on all walls to ensure the flow tangency. The initial conditions are as follows:
\begin{subequations}
    \begin{align}
        T &= 35 \; \text{K}, \; M_x = 10, \; M_y = -3, \; p = 0.1 \; \text{bar} &\forall x \leq 0 \\
        T &= 35 \; \text{K}, \; M_x = 10, \; M_y = 2, \; p = 0.1 \; \text{bar} &\forall x \geq 0
    \end{align}
\end{subequations}
The CFL number is fixed at 0.20, and the flow is allowed to evolve for $ 9.5 \times 10^{-4} $ time units. Figure \ref{2D_figure5} shows the logarithmic pressure contours computed using the second-order accurate modified PVU Boltzmann scheme. This numerical result reinforces the theoretical proof that the second-order accurate modified PVU Boltzmann scheme preserves the positivity of pressure and density.

\subsubsection{Shock diffraction}
\label{subsubsec 9.3.4}

\begin{figure}[H]
   \centering
   \includegraphics[scale=0.3]{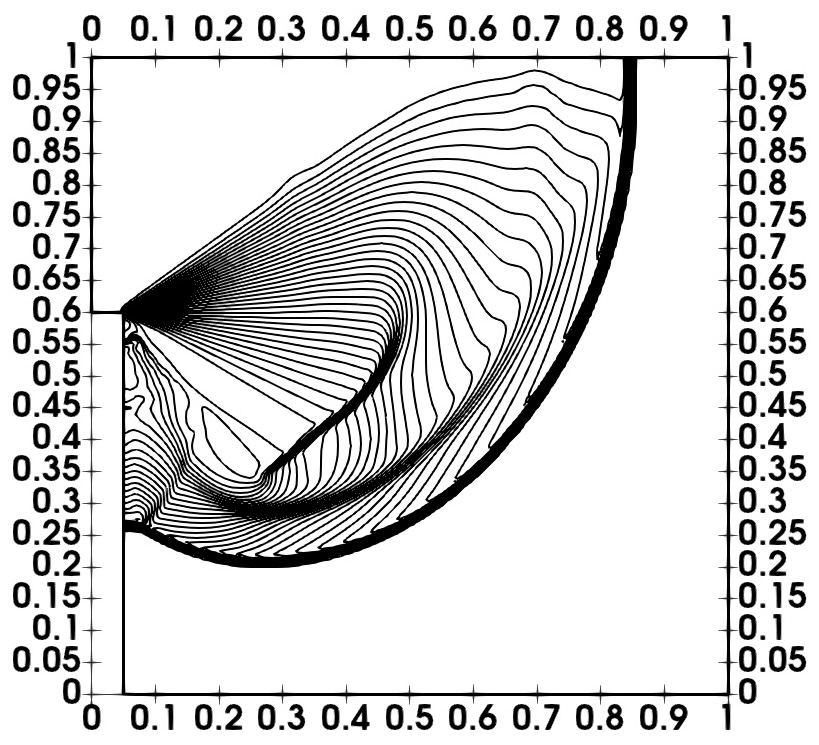}
   \caption{Shock diffraction; density contours (0.01:0.177:7.1)}
   \label{2D_figure4}
\end{figure}

In this test case, a normal shock of Mach 5.09 diffracts around a $ 90^{\circ} $ corner. The computational domain is a square of $ 1 $ length units, which is divided into $ 200 \times 200 $ Cartesian cells. The corner is located at $ x = 0.05 $ and $ y = 0.45 $ length units. The initial location of the shock is at the corner. The domain to the right of the shock is initialized with the pre-shock conditions of $ \rho = 1.4, \, u_1 = 0.0, \, u_2 = 0.0, \, p = 1.0 $. The domain to the left of the shock is initialized using the post-shock conditions. The boundary conditions used are the supersonic inlet boundary condition at the left boundary, the kinetic characteristic boundary condition at the walls, and the zero gradient boundary condition at the right boundary of the domain. The solution evolves until $ 0.1561 $ time units. The second order results corresponding to the modified PVU Boltzmann scheme are shown in Figure \ref{2D_figure4}. As seen in the figure, the numerical results are free from any entropy violating solutions.

\subsubsection{Supersonic flow over a forward-facing step}
\label{subsubsubsec 9.3.5}

\begin{figure}[ht]
   \centering
   \includegraphics[scale=0.27]{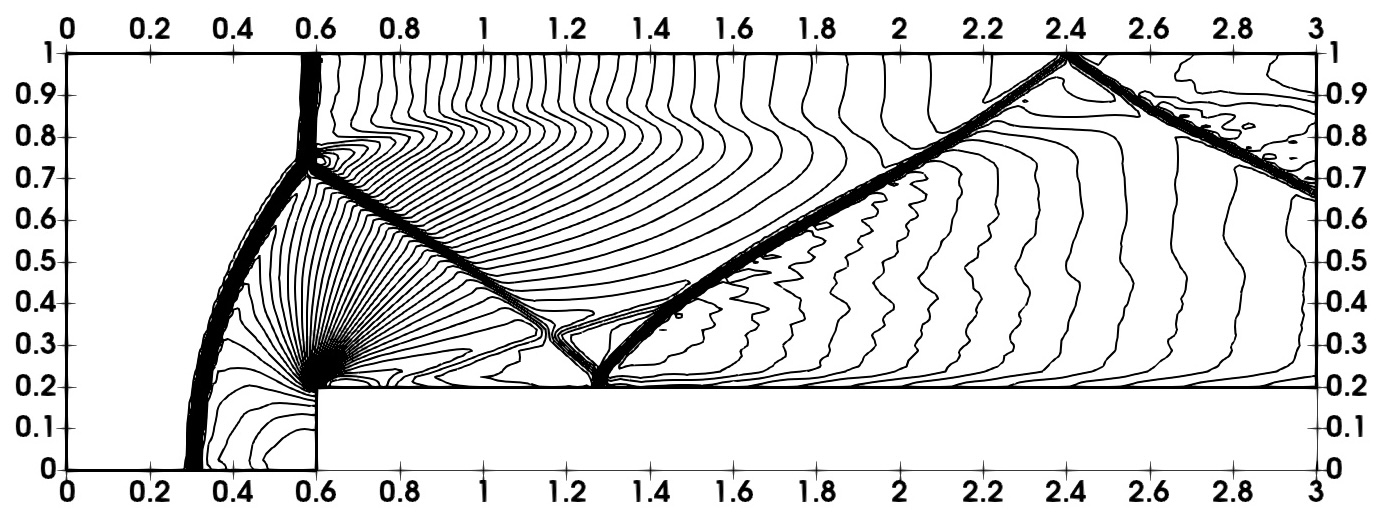}
   \caption{Supersonic flow over a forward-facing step; density contours (0.15:0.13:6.30)}
   \label{2D_figure6}
\end{figure}
This problem consists of a supersonic flow of Mach 3 in a wind tunnel with a step. The wind tunnel is $ 1 $ length units wide and $ 3 $ length units long. The step is $ 0.2 $ length units in height and is placed at $ 0.6 $ length units from the left end of the tunnel. The computational domain is divided into $ 240 \times 80 $ Cartesian cells. The supersonic inlet boundary condition is assumed at the left end, the kinetic characteristic boundary condition is used at the top and bottom walls, and the zero gradient boundary condition is assumed at the right boundary of the tunnel.  The initial conditions corresponding to the free stream are assumed everywhere in the cell. The flow evolves until $ t = 4 $ time units. The second order results for this numerical test computed using the proposed scheme are shown in Figure \ref{2D_figure6}.

\subsubsection{Double mach reflection}
\label{subsubsec 9.3.6}

\begin{figure}[h]
   \centering
   \includegraphics[scale=0.24]{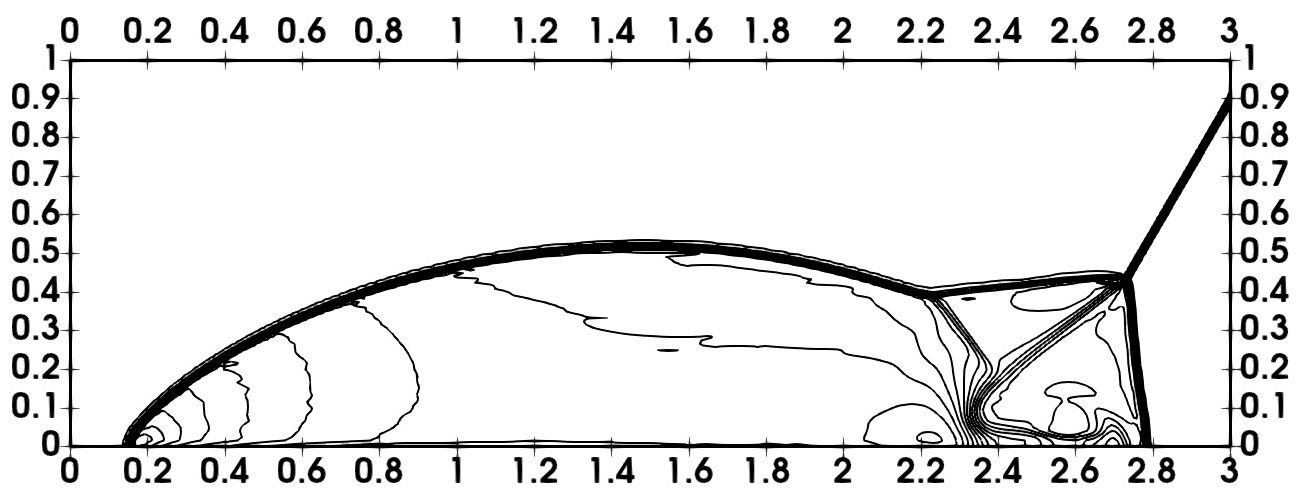}
   \caption{Double mach reflection; density contours (2.0:0.65:21.5)}
   \label{2D_figure7}
\end{figure}

This test case is designed to assess the performance of a numerical scheme with respect to another shock instability known as a kinked mach stem. The problem is set up by keeping the domain rectangular and forcing an oblique shock through it \cite{Woodward_Colella}. The computational domain is four units long and one unit wide. It is divided into $ 480 \times 160 $ Cartesian cells. An oblique Mach 10 shock that initially has an angle of $ 60^{\circ} $ with the bottom wall at $ x = 1/6 $ is propagated in a channel. The domain to the right of the shock is initialized with pre-shock values of $ \rho = 1.4, u_1 = 0.0, u_2 = 0.0, p = 1.0 $, and to the left is initialized using post-shock values. The inlet boundary condition is set to post-shock values, and the outlet boundary is set to zero gradient. The top boundary condition is set to simulate the actual shock movement. The bottom wall up to $ x = 1/6 $ is set to post-shock values and the kinetic characteristic boundary condition is used for the remaining wall. The flow evolves until $ 0.2 $ time units. Some low diffusive numerical schemes are plagued by the unphysical kinked mach stem.  Some low diffusive numerical schemes are plagued by the unphysical kinked mach stem. As seen in Figure \ref{2D_figure7}, the modified PVU Boltzmann scheme does not produce such an unphysical kink in the Mach stem.

\subsubsection{Hypersonic flow over a half cylinder}
\label{subsubsec 9.3.7}

\begin{figure}[h]
   \centering
   \includegraphics[scale=0.3]{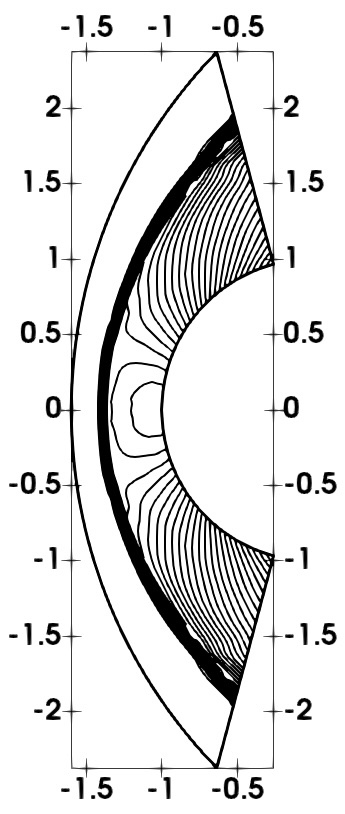}
   \caption{Hypersonic flow over a half cylinder; density contours (1.3:0.25:8.9)}
   \label{2D_figure8}
\end{figure}
This numerical test consists of a hypersonic flow of Mach 20 over a half cylinder. This test is useful in assessing the performance of the scheme with respect to carbuncle instability \cite{Quirk}. We use the algorithm given in \cite{Huang_Wu} to generate the mesh for this test case. The computational domain is divided into $ 160 \times 20 $ cells. The domain is initialized with the free-flow values of $ \rho = 1.4, u_1 = 20.0, u_2 = 0.0, p = 1.0 $. The supersonic inlet boundary condition is used at the left boundary of the domain, the kinetic characteristic boundary condition is used at the cylinder surface, and the zero gradient boundary condition is used at the outlet of the domain. In the steady regime the flow field exhibits a detached bow shock standing upstream of the half-cylinder. As seen in Figure \ref{2D_figure8}, the numerical results for this test case computed using the second-order modified PVU Boltzmann scheme are free from the Carbuncle instability.


\section{Conclusions}
\label{sec 10}

This paper has presented the formulation and positivity analysis of the peculiar velocity based Boltzmann scheme derived using the compactly supported hat function. This development lays the groundwork for the derivation of the Boltzmann scheme to capture contact discontinuities, introduced in \cite{Singh_Rao_2026}, using the hat function. Next, we present the formulation and positivity analysis of the contact discontinuity capturing Boltzmann scheme based on the hat function. Building on this foundation, we inherit ideas from \cite{Zhang_Shu} to extend the proposed kinetic scheme to second-order accuracy in space and time, while preserving the positivity property of the underlying first-order formulation. The associated positivity analysis yields a CFL-like condition for selecting the time step in second-order simulations. Finally, we demonstrate the robustness and positivity-preserving nature of the modified PVU Boltzmann scheme through a series of one- and two-dimensional compressible flow problems.


\section*{CRediT author statement}

\textbf{Balwinder Singh:} Methodology, Conceptualization, Investigation, Software, Validation, Formal analysis, Writing - original draft.

\vspace{1 em}
\textbf{S. V. Raghurama Rao:} Conceptualization, Investigation, Formal Analysis, Supervision, Writing - Review \& Editing.

\section*{Declaration of competing interest}

We confirm that neither monetary interests nor individual relationships involving persons or institutions that could influence the material presented in this study.

\section*{Data availability}

This article does not rely on any datasets for the research it presents.

\section*{Acknowledgments}

This research did not receive specific grants from funding agencies in the public, commercial or non-profit sectors.


\appendix

\section{Details of the 1D numerical flux}
\noindent
\textbf{PVU:}
\begin{equation}
    \hat{\textbf{f}}_{j + \frac{1}{2}}^{(p)} = \frac{1}{2} \left( \textbf{G}_j^{(p)} + \textbf{G}_{j+1}^{(p)} \right) - \textbf{d}_{j + \frac{1}{2}}^{(p)}
\end{equation}
where $ \textbf{d}_{j + \frac{1}{2}}^{(p)} $ denotes the diffusive pressure flux, which reads
\renewcommand{\arraystretch}{2.5}
\begin{equation}
    \textbf{d}_{j + \frac{1}{2}}^{(p)} = \begin{bmatrix}
        \dfrac{1}{6} \Delta_{j + \frac{1}{2}} \left( \rho \, \sqrt{6 R T} \right) \\
        \dfrac{1}{6} \Delta_{j + \frac{1}{2}} \left( \rho \, u \, \sqrt{6 R T} \right) \\ 
        \dfrac{1}{12} \Delta_{j + \frac{1}{2}} \left( \rho \, u^2 \, \sqrt{6 R T} \right) + \dfrac{1}{10} \left( \dfrac{\gamma + 1}{\gamma - 1} \right) \Delta_{j + \frac{1}{2}} \left( p \sqrt{6 R T} \right)
    \end{bmatrix}.
\end{equation}

\vspace{1\baselineskip}
\noindent
\textbf{PVU+:}
The diffusive flux pertaining to numerical pressure flux of the modified PVU Boltzmann scheme reads
\begin{equation}
    \textbf{d}_{j + \frac{1}{2}}^{(p)} = \begin{bmatrix}
        \Delta_{j + \frac{1}{2}} \left( \dfrac{\rho \, u^2}{4 \, \sqrt{3 R T}} \left( 1 - \exp{\left( \dfrac{-3 R T}{u^2} \right)} \right) \right) \\
        \Delta_{j + \frac{1}{2}} \left( \dfrac{\rho \, u^3}{4 \, \sqrt{3 R T}} \left( 1 - \exp{\left( \dfrac{-3 R T}{u^2} \right)} \right) \right) \\
        \begin{array}{cc}
             \Delta_{j + \frac{1}{2}} \left( \dfrac{\rho \, u^4}{4 \, \sqrt{3 R T}} \left( 1 - \exp{\left( \dfrac{-3 R T}{u^2} \right)} \left( 1 + \dfrac{3 R T}{u^2} \right) \right) \right) \\
             + \, \Delta_{j + \frac{1}{2}} \left( \dfrac{u^2}{2} \, G_1^+ \right) + \Delta_{j + \frac{1}{2}} \left( \dfrac{\left( 3 - \gamma \right) R T}{2 \left( \gamma - 1 \right)} G_1^+ \right)
        \end{array}
    \end{bmatrix},
\end{equation}
where
\begin{equation}
    G_1^+ = \dfrac{\rho \, u^2}{4 \, \sqrt{3 R T}} \left( 1 - \exp{\left( \dfrac{-3 R T}{u^2} \right)} \right).
\end{equation}

\section{Details of the 2D numerical flux}
\noindent
\textbf{PVU:}
\begin{equation}
    \hat{\textbf{f}}_{e}^{(p)} = \frac{1}{2} \left( \textbf{G}_{n, K}^{(p)} + \textbf{G}_{n, K_e}^{(p)} \right) - \textbf{d}_e^{(p)},
\end{equation}
where $ \textbf{d}_{e}^{(p)} $ denotes the diffusive pressure flux, which reads
\renewcommand{\arraystretch}{2.5}
\begin{equation}
    \textbf{d}_e^{(p)} = \begin{bmatrix}
        \dfrac{1}{6} \Delta_e \left( \rho \, \sqrt{6 R T} \right) \\
        \dfrac{1}{6} \Delta_e \left( \rho \, u_n \, \sqrt{6 R T} \right) \\
        \dfrac{1}{6} \Delta_e \left( \rho \, u_t \, \sqrt{6 R T} \right) \\
        \dfrac{1}{12} \Delta_e \left( \rho \, | \boldsymbol{u} |^2 \, \sqrt{6 R T} \right) + \dfrac{1}{10} \left( \dfrac{\gamma + 1}{\gamma - 1} \right) \Delta_e \left( p \sqrt{6 R T} \right)
    \end{bmatrix}.
\end{equation}

\vspace{1\baselineskip}
\noindent
\textbf{PVU+:}
\begin{equation}
    \textbf{d}_{e}^{(p)} = \begin{bmatrix}
        \Delta_e \left( \dfrac{\rho \, u_n^2}{4 \, \sqrt{3 R T}} \left( 1 - \exp{\left( \dfrac{-3 R T}{u_n^2} \right)} \right) \right) \\
        \Delta_e \left( \dfrac{\rho \, u_n^3}{4 \, \sqrt{3 R T}} \left( 1 - \exp{\left( \dfrac{-3 R T}{u_n^2} \right)} \right) \right) \\
        \Delta_e \left( \dfrac{\rho \, u_n^2 \, u_t}{4 \, \sqrt{3 R T}} \left( 1 - \exp{\left( \dfrac{-3 R T}{u_n^2} \right)} \right) \right) \\
        \begin{array}{cc}
             \Delta_e \left( \dfrac{\rho \, u_n^4}{4 \, \sqrt{3 R T}} \left( 1 - \exp{\left( \dfrac{-3 R T}{u_n^2} \right)} \left( 1 + \dfrac{3 R T}{u_n^2} \right) \right) \right) \\
             + \, \Delta_e \left( \dfrac{| \boldsymbol{u} |^2}{2} \, G_1^+ \right) + \Delta_e \left( \dfrac{\left( 3 - \gamma \right) R T}{2 \left( \gamma - 1 \right)} G_1^+ \right)
        \end{array}
    \end{bmatrix},
\end{equation}
where
\begin{equation}
    G_1^+ = \dfrac{\rho \, u_n^2}{4 \, \sqrt{3 R T}} \left( 1 - \exp{\left( \dfrac{-3 R T}{u_n^2} \right)} \right).
\end{equation}


\end{document}